\documentclass[5p,twocolumn,number,fleqn,A4]{elsarticle}

\usepackage{color} 
\usepackage{graphicx} 
\usepackage{latexsym,amssymb,amsmath} 
\usepackage{mathrsfs} 
\usepackage{multicol}
\usepackage{fixltx2e} 
\usepackage{cancel} 

\newtheorem{thm}{Theorem}
\newtheorem{lem}{Lemma}
\newtheorem{defn}{Definition}
\newtheorem{exam}{Example}
\newtheorem{coro}{Corollary}
\newproof{pf}{Proof}

\allowdisplaybreaks[4]

\def\abs#1{\left\vert #1 \right\vert}
\def\allpoly{\mbox{$\re\langle X \rangle$}}

\def\allpolyx0degn{\mbox{$P_n$}}

\def\allseries{\mbox{$\re\langle\langle X \rangle\rangle$}}

\def\allseriesXcmd{\mbox{$\re^{m_d}\langle\langle X_c \rangle\rangle$}}
\def\allseriesXcmdLC{\mbox{$\re^{m_d}_{LC}\langle\langle X_c \rangle\rangle$}}

\def\allseriesdeltam{\mbox{$\re^m\langle\langle X_\delta \rangle\rangle$}}

\def\allseriesdeltamLC{\mbox{$\re^m_{LC}\langle\langle X_\delta \rangle\rangle$}}
\def\allseriesXdmc{\mbox{$\re^{m_c}\langle\langle X_d \rangle\rangle$}}
\def\allseriesXdmcLC{\mbox{$\re^{m_c}_{LC}\langle\langle X_d \rangle\rangle$}}

\def\allseriesell{\mbox{$\re^{\ell} \langle\langle X \rangle\rangle$}}

\def\allseriesm{\mbox{$\re^m\langle\langle X \rangle\rangle$}}

\def\allseriesmGC{\mbox{$\re^{m}_{GC}\langle\langle X \rangle\rangle$}}
\def\allseriesmLC{\mbox{$\re^{m}_{LC}\langle\langle X \rangle\rangle$}}

\def\allseriesellLC{\mbox{$\re^{\ell}_{LC}\langle\langle X \rangle\rangle$}}
\def\allseriesellGC{\mbox{$\re^{\ell}_{GC}\langle\langle X \rangle\rangle$}}

\def\allseriesX1{\mbox{$\re [[ X_1 ]]$}}
\newcommand{\comment}[1]{} 
\def\conc{{\rm cat}}

\def\dim{{\rm dim}}

\def\Endallseries{{\rm End}(\allseries)}
\def\eqref#1{(\ref{#1})} 

\def\Fliessdelta{\mathscr{F}_{\delta}}

\def\id{{\rm id}}

\def\Lpm{L_{\mathfrak{p}}^m}

\def\Lpme{L^m_{\mathfrak{p},e}}

\def\mbf#1{\hbox{\mathversion{bold}$#1$}} 
\def\modcomp{\:\tilde{\circ}\,} 

\def\norm#1{\Vert#1\Vert}

\def\re{{\mathbb R}} 

\def\sameau{\rule[0.017in]{0.2in}{0.012in}}
\def\sh{{\rm sh}}
\def\shuffle{{\scriptscriptstyle \;\sqcup \hspace*{-0.05cm}\sqcup\;}}

\def\spanset{{\rm span}}

\def\gsc{c} 
\def\gsd{d} 

\def\begce{\begin{center}}
\def\endce{\end{center}}
\def\begar{\begin{array}}
\def\endar{\end{array}}
\def\begeq{\begin{equation}}
\def\endeq{\end{equation}}
\def\begdi{\begin{displaymath}}
\def\enddi{\end{displaymath}}
\def\begdis{\begin{eqnarray*}}
\def\enddis{\end{eqnarray*}}
\def\begeqa{\begin{eqnarray}}
\def\endeqa{\end{eqnarray}}
\def\begdes{\begin{description}}
\def\enddes{\end{description}}
\def\begit{\begin{itemize}}
\def\endit{\end{itemize}}
\def\begen{\begin{enumerate}}
\def\enden{\end{enumerate}}
\def\beglar{\left[\begin{array}}
\def\endrar{\end{array}\right]}
\def\begle{\begin{lem}}
\def\endle{\end{lem}}
\def\begde{\begin{defn}}
\def\endde{\end{defn}}
\def\begth{\begin{thm}}
\def\endth{\end{thm}}
\def\begco{\begin{coro}}
\def\endco{\end{coro}}
\def\begprop{\begin{proposition}}
\def\endprop{\end{proposition}}
\def\begex{\begin{exam}}
\def\endex{\end{exam}}
\def\begexer{\begin{exercise}}
\def\endexer{\end{exercise}}

\def\begres{\noindent{\bf Remarks}:\begin{enumerate}}
\def\endres{\end{enumerate} \par}
\def\begpr{\begin{pf}}
\def\endpr{\end{pf}}
\def\begtab{\begin{tabular}}
\def\endtab{\end{tabular}}
\def\rref#1{(\ref{#1})}

\newcommand{
\markblue}[1]{{\textcolor[rgb]{0,0,0}{#1}}} 

\begin{document}

\begin{frontmatter}

\title{Fa\`{a} di Bruno Hopf Algebra of the Output \\
Feedback Group for Multivariable Fliess Operators}

\author[icmat,odu]{W.~Steven Gray\corref{cor}}
\ead{sgray@odu.edu}

\author[gmu]{Luis A.~Duffaut Espinosa}
\ead{lduffaut@gmu.edu}

\author[icmat]{Kurusch Ebrahimi-Fard}
\ead{kurusch@icmat.es}

\cortext[cor]{Corresponding author. Tel.:~+34-91-29-99-719}

\address[icmat]{Instituto de Ciencias Matem\'{a}ticas,
Consejo Superior de Investigaciones Cient\'{\i}ficas,
C/ Nicol\'{a}s Cabrera, no.~13-15,
28049 Madrid, Spain}

\address[odu]{On leave from Old Dominion University, Norfolk, Virginia 23529, USA}

\address[gmu]{Department of Electrical and Computer Engineering, George Mason University,
Fairfax, Virginia 22030, USA}

\begin{abstract}
Given two nonlinear input-output systems written in terms
of Chen-Fliess functional expansions, i.e., Fliess operators, it is known
that the feedback interconnected system
is always well defined and in the same class. An explicit
formula for the generating series of a single-input, single-output
closed-loop system was
provided by the first two authors in earlier work via
Hopf algebra methods.
This paper is a sequel. It has four main innovations.
First, the full multivariable extension of the theory is
presented.
Next, a major simplification of the basic setup is introduced
using a new type of grading that has recently appeared in the literature.
This grading also facilitates a fully recursive
algorithm to compute the antipode of the Hopf algebra of the output feedback group,
and thus, the corresponding
feedback product can be computed much more efficiently.
The final innovation is an improved convergence analysis of
the antipode operation, namely, the radius of convergence of
the antipode is computed.
\end{abstract}


\end{frontmatter}

\section{Introduction}

Given two nonlinear input-output systems written in terms
of Chen-Fliess functional expansions \cite{Fliess_81}, i.e., Fliess operators, it was shown in
\cite{Gray-Li_05,Gray-Wang_08} that the feedback interconnected system
is always well defined and in the same class. An explicit
formula for the generating series of a single-input, single-output (SISO)
closed-loop system was
later provided in \cite{Gray-Duffaut_Espinosa_SCL11} using
Hopf algebra methods. In particular, the so called {\em feedback product}
of the two generating series for the component systems can be computed
in terms of the antipode of a Fa\`{a} di Bruno type Hopf algebra.
This antipode was described in terms of a sequence of polynomials of
increasing degree. While explicit, this somewhat brute force formula is not ideal for software
implementation \cite{Gray-et-al_MTNS14}. Nevertheless, this antipode can be used to
provide a tractable formula for nonlinear system inversion from a purely
input-output point of view, i.e., no state space model is
required \cite{Gray-et-al_AUT}.

This paper is a sequel to \cite{Gray-Duffaut_Espinosa_SCL11}. It has four main innovations.
First, the full multivariable extension of the theory in \cite{Gray-Duffaut_Espinosa_SCL11} is
presented, which makes it more relevant to practical control problems.
The second innovation is more technical, but it greatly simplifies the basic setup.
Specifically, it was shown recently in \cite{Foissy_13} that the Hopf algebra for the SISO
output feedback
group is {\em connected} under
a grading that is distinct from the one described in \cite{Gray-Duffaut_Espinosa_SCL11}. This
important observation implies that the bialgebra presented in the original paper is {\em automatically} a
Hopf algebra, and therefore, much of the technical analysis concerning the existence of
the antipode can now be omitted. So here the method in \cite{Foissy_13} is extended to the
multivariable case and applied throughout.
The third innovation is related to the existence
of this new grading. Namely,
the {\em partially} recursive formula for the antipode of any connected graded Hopf algebra
in \cite{Figueroa-Gracia-Bondia_05} is exploited here to
produce a {\em fully} recursive antipode
algorithm for the Hopf algebra of the output feedback group. This in turn allows one to compute the feedback product
much more efficiently. The approach involves carefully
combining results from
\cite{Figueroa-Gracia-Bondia_05}, \cite{Foissy_13} and \cite{Reutenauer_93}.
The SISO version of this algorithm was presented in \cite{Gray-et-al_MTNS14} and
compared against other existing methods.  In a Mathematica implementation, this new algorithm provided
an order of magnitude reduction in execution times. For the multivariable case, such gains are likely to be even larger, but this analysis
is beyond the scope of this paper. The final innovation is an improved convergence analysis of
the antipode operation, specifically, the radius of convergence of the antipode is computed
using techniques presented in \cite{Thitsa-Gray_12}.
In \cite{Gray-Duffaut_Espinosa_SCL11} it was only shown that this radius of convergence
is positive.

The paper is organized as follows. In the next section, some mathematical preliminaries and background
are \markblue{briefly summarized concerning Hopf algebras and the interconnection of Fliess
operators.} In Section~\ref{sec:FdB-Hopf-algebra}, the Hopf algebra of the multivariable
output feedback group is presented, including the recursive algorithm for the antipode
and the radius of convergence for this operation.
In the subsequent section, these results are used to define the multivariable
feedback product, and the corresponding convergence analysis is presented.
The theory is demonstrated on a simple steering example.
The conclusions are given in the final section.

\section{Preliminaries}
 \label{sec:preliminaries}

A finite nonempty set of noncommuting symbols $X=\{ x_0,x_1,\ldots,x_m\}$ is
called an {\em alphabet}. Each element of $X$ is called a {\em
letter}, and any finite sequence of letters from $X$,
$\eta=x_{i_1}\cdots x_{i_k}$, is called a {\em word} over $X$. The
{\em length} of $\eta$, $\abs{\eta}$, is the number of letters in
$\eta$.
Let $\abs{\eta}_{x_i}$ denote the number of times the letter $x_i\in X$
appears in the word $\eta$.
The set of all words including the empty word, $\emptyset$,
is designated by $X^\ast$. It forms a monoid under catenation.
Any mapping $c:X^\ast\rightarrow
\re^\ell$ is called a {\em formal power series}. The value of $c$ at
$\eta\in X^\ast$ is written as $(c,\eta)$ and called the {\em coefficient} of
$\eta$ in $c$.
Typically, $c$ is
represented as the formal sum $c=\sum_{\eta\in X^\ast}(c,\eta)\eta.$
The collection of all formal power series over $X$ is denoted by
$\allseriesell$.
It forms an associative $\re$-algebra under
the catenation product and an associative and
commutative $\re$-algebra under the {\em shuffle product}, that is,
the bilinear product defined in terms of the shuffle product of two words
\begdi
(x_i\eta)\shuffle(x_j\xi)=x_i(\eta\shuffle(x_j\xi))+x_j((x_i\eta)\shuffle \xi),
\enddi
where $x_i,x_j\in X$, $\eta,\xi\in X^\ast$ and with
$\eta\shuffle\emptyset=\emptyset\shuffle\eta=\eta$ \cite{Fliess_81,Reutenauer_93}.
Its restriction to polynomials over $X$ is
\begdi
\sh : \allpoly\otimes\;\allpoly \rightarrow \allpoly
    : p\otimes q \mapsto p\shuffle q.
\enddi
The corresponding adjoint map $\sh^{\ast}$ is
the unique
$\re$-linear map of the form $\allpoly\rightarrow \allpoly\otimes\allpoly$
which satisfies the identity
\begdi
(\sh(p\otimes q),r)=(p\otimes q,\sh^{\ast}(r)) \label{eq:sh}
\enddi
for all $p,q,r\in\allpoly$.
The following theorem states an important duality.

\begth \label{th:shuffle_morphs} \cite{Reutenauer_93}
The adjoint map $\sh^{\ast}$ is an $\re$-algebra
morphism for the catenation product $\conc:p\otimes q\mapsto pq$. That is,
\begdi
\sh^{\ast}(pq)=\sh^{\ast}(p)\,\sh^{\ast}(q)
\label{eq:shstar-duality}
\enddi
for all $p,q\in \allpoly$ with $\sh^\ast(\mbf{1})=\mbf{1}\otimes \mbf{1}$.\footnote{Here $\mbf{1}$ is the unit polynomial $1\emptyset$.}
In particular, for $x_i\in X$ and $\eta\in X^{\ast}$
\begdi
\sh^{\ast}(x_i\eta)=(x_i\otimes \mbf{1} + \mbf{1}\otimes x_i)\sh^{\ast}(\eta). \label{eq:shstar-on-Xstar} \\[0.1in]
\enddi
\endth

\subsection{\markblue{Hopf Algebras}}

\markblue{
In this section, a few basic facts and tools concerning Hopf algebras are summarized.
The reader is referred to \cite{Figueroa-Gracia-Bondia_05,Hochschild_81,Sweedler_69} for more complete treatments.
}

\markblue{
A {\it{coalgebra}} over $\re$ consists of a triple $(C,\Delta,\varepsilon)$. The coproduct $\Delta: C \to C \otimes C$ is coassociative, that is,
$(\id \otimes \Delta)\circ \Delta=(\Delta \otimes \id)\circ \Delta$, and $\varepsilon: C \to \re$ denotes the counit map. A {\it{bialgebra}} $B$ is both a unital algebra and a coalgebra together with compatibility relations, such as both the algebra product, $m(x,y)=xy$, and unit map, $\mathsf{e}: \re \to B$, are coalgebra morphisms. This provides, for example, that
$\Delta(xy)=\Delta(x)\Delta(y)$.
The unit of $B$ is denoted by $\mathbf{1} = \mathsf{e}(1)$. A bialgebra is called {\it{graded}} if there are $\re$-vector subspaces $B_n$,
$n \geq 0$ such that $B= \bigoplus_{n \geq 0} B_n$ with $m(B_k \otimes B_l) \subseteq B_{k+l}$
and $\Delta B_n \subseteq \bigoplus_{k+l=n} B_k\otimes B_l.$ Elements $x \in B_n$ are given a degree
$\deg(x)=n$. Moreover, $B$ is called {\it{connected}} if $B_0 = \re\mathbf{1}$. Define $B_+=\bigoplus_{n > 0} B_n$. For any $x \in B_n$ the coproduct is of the form
\begdi
	\Delta(x) = x \otimes \mathbf{1} + \mathbf{1} \otimes x + \Delta'(x) \in \bigoplus_{k+l=n} B_k \otimes B_l,
\enddi
where $\Delta'(x) := \Delta(x) -  x \otimes \mathbf{1} - \mathbf{1} \otimes x \in B_+\otimes B_+$ is the {\em reduced} coproduct.
}

\markblue{
Suppose $A$ is an $\re$-algebra with product $m_A$ and unit $\mathsf{e}_A$, e.g., $A=\re$ or $A=B$. The vector space  $L(B, A)$ of linear maps from the bialgebra $B$ to $A$ together with the convolution product $\Phi \star \Psi := m_{A} \circ (\Phi \otimes \Psi) \circ \Delta : B \to A$, where $\Phi,\Psi \in L(B,A)$, is an associative algebra with unit $\iota := \mathsf{e}_{A} \circ \varepsilon$.
A {\it{Hopf algebra}} $H$ is a bialgebra together with a particular $\re$-linear map called an {\it{antipode}} $S: H \to H$
which satisfies the Hopf algebra axioms and has the property that $S(xy)=S(y)S(x)$. When $A=H$, the antipode $S \in L(H,H)$ is the inverse of the identity map with respect to the convolution product, that is,
\begin{equation*}
    S  \star\id = \id \star S := m \circ (S \otimes \id) \circ \Delta = \mathsf{e} \circ \varepsilon.
\end{equation*}
A connected graded bialgebra $H=\bigoplus_{n \ge 0} H_n$ is {\em always} a connected graded Hopf algebra.
}

\markblue{
Suppose $A$ is a commutative unital algebra. The subset $g_0 \subset  L( H, A)$ of linear maps $\alpha$ satisfying $\alpha(\mathbf{1})=0$ forms a Lie algebra in $ L( H, A)$. The exponential $ \exp^\star(\alpha) = \sum_{j\ge 0} \frac{1}{j!}\alpha^{\star j}$ is well defined and gives a bijection from $g_0$ onto the group $G_0 = \iota + g_0$ of linear maps $\gamma$ satisfying $\gamma(\mathbf{1})=1_{A}$.
A map $\Phi\in L( H, A)$ is called a {\it{character}} if $\Phi(\mathbf{1})=1_{ A}$ and $\Phi(xy) = \Phi(x)\Phi(y)$ for all $x,y\in H$.
The set of characters is denoted by $G_{ A} \subset G_0$.
The neutral element $\iota:=\mathsf{e}_{ A}\circ \varepsilon$  in $G_{ A}$ is given by $\iota(\mathbf{1})=1_{ A}$ and $\iota(x) = 0$ for $x \in \mathsf{Ker}(\varepsilon)=H_+$. The inverse of $\Phi \in G_{ A}$ is given by
\begeq \label{eq:convinv-via-antipode}
\Phi^{\star -1} = \Phi \circ S.
\endeq}%

\vspace{-0.2in}

\markblue{Given an arbitrary group $G$, the set of real-valued functions defined on $G$ is a commutative unital algebra. There is a subalgebra of functions known as the
{\em representative functions}, $R(G)$, which can be endowed with a Hopf algebra, $H$. In this case, there is a group isomorphism relating $G$ to the convolution
group $G_A$, say, $\Phi:G\rightarrow G_A:g\mapsto \Phi_g$.
A {\em coordinate map} is any $a:H\rightarrow \re$ satisfying
\begeq \label{eq:general-coordinate-maps}
(\Phi_{g_1}\star \Phi_{g_2})(a)=a(g_1g_2),\;\;\forall g_i\in G.
\endeq
In some sense, the coordinates maps are the generators of $H$, though they can not always be easily identified in general.
}%

\vspace{-0.1in}

\markblue{
\begex \cite{Reutenauer_93}
$(\allpoly,\conc,\mathsf{e},\sh^\ast,\varepsilon,S)$
is a Hopf algebra, where
$\mathsf{e}:\re\rightarrow \allpoly : k\mapsto k {\mbf 1}$,
$\varepsilon:\allpoly\rightarrow\re: p\mapsto (p,\emptyset)$,
\begdi
f\star g:\allpoly\rightarrow\allpoly
:p \mapsto \sum_{\eta,\xi\in X^\ast} (p,\eta\shuffle\xi)\,f(\eta)g(\xi),
\enddi
for all $f,g\in L(\allpoly,\allpoly)$,
and
$
S(x_{i_1}x_{i_2}\cdots x_{i_k})=(-1)^k x_{i_k}x_{i_{k-1}}\cdots x_{i_1}.
$
\endex
}

\subsection{Fliess Operators and Their Interconnections}

One can formally associate with any series $c\in\allseriesell$ a causal
$m$-input, $\ell$-output operator, $F_c$, in the following manner.
Let $\mathfrak{p}\ge 1$ and $t_0 < t_1$ be given. For a Lebesgue measurable
function $u: [t_0,t_1] \rightarrow\re^m$, define
$\norm{u}_{\mathfrak{p}}=\max\{\norm{u_i}_{\mathfrak{p}}: \ 1\le
i\le m\}$, where $\norm{u_i}_{\mathfrak{p}}$ is the usual
$L_{\mathfrak{p}}$-norm for a measurable real-valued function,
$u_i$, defined on $[t_0,t_1]$.  Let $L^m_{\mathfrak{p}}[t_0,t_1]$
denote the set of all measurable functions defined on $[t_0,t_1]$
having a finite $\norm{\cdot}_{\mathfrak{p}}$ norm and
$B_{\mathfrak{p}}^m(R)[t_0,t_1]:=\{u\in
L_{\mathfrak{p}}^m[t_0,t_1]:\norm{u}_{\mathfrak{p}}\leq R\}$. Assume $C[t_0,t_1]$
is the subset of continuous functions in $L_{1}^m[t_0,t_1]$. Define
inductively for each $\eta\in X^{\ast}$ the map $E_\eta:
L_1^m[t_0, t_1]\rightarrow C[t_0, t_1]$ by setting
$E_\emptyset[u]=1$ and letting
\[E_{x_i\bar{\eta}}[u](t,t_0) =
\int_{t_0}^tu_{i}(\tau)E_{\bar{\eta}}[u](\tau,t_0)\,d\tau, \] where
$x_i\in X$, $\bar{\eta}\in X^{\ast}$, and $u_0=1$. The
input-output operator corresponding to $c$ is the {\em Fliess operator}
\begeq
F_c[u](t) =
\sum_{\eta\in X^{\ast}} (c,\eta)\,E_\eta[u](t,t_0) \label{eq:Fliess-operator-defined}
\endeq
\cite{Fliess_81,Fliess_83}.
If there exist real numbers $K_c,M_c>0$ such that
\begeq
\abs{(c,\eta)}\le K_c M_c^{|\eta|}|\eta|!,\;\; \forall\eta\in X^{\ast},
\label{eq:local-convergence-growth-bound}
\endeq
then $F_c$ constitutes a well defined mapping from
$B_{\mathfrak p}^m(R)[t_0,$ $t_0+T]$ into $B_{\mathfrak
q}^{\ell}(S)[t_0, \, t_0+T]$ for sufficiently small $R,T >0$,
where the numbers $\mathfrak{p},\mathfrak{q}\in[1,\infty]$ are
conjugate exponents, i.e., $1/\mathfrak{p}+1/\mathfrak{q}=1$
\cite{Gray-Wang_SCL02}.
(Here, $\abs{z}:=\max_i \abs{z_i}$ when $z\in\re^\ell$.) The set of all such
{\em locally convergent} series is denoted by $\allseriesellLC$.
In particular, when $\mathfrak{p}=1$,
the series \rref{eq:Fliess-operator-defined}
converges absolutely and uniformly if
$
\max\{R,T\}< 1/M_c(m+1)
$
\cite{Duffaut-Espinosa_09,Duffaut-Espinosa-et-al_CDC09}.
It is important in applications to identify the {\em smallest} possible
geometric growth constant, $M_c$, in order to avoid over restricting the domain of
$F_c$.
So let $\pi:\allseriesellLC\rightarrow\re^+\cup \; \{0\}$
take each series $c$
to the infimum of all $M_c$ satisfying
\rref{eq:local-convergence-growth-bound}. Therefore,
$\allseriesellLC$ can be partitioned into equivalence classes, and
the number $1/M_c(m+1)$ will be referred to as the {\em radius of convergence} for
the class
$\pi^{-1}(M_c)$. This is in contrast to the usual situation where a radius
of convergence is assigned to individual series.
When $c$ satisfies the more stringent growth condition
\begeq
\abs{(c,\eta)}\le K_c M_c^{|\eta|},\;\; \forall\eta\in X^{\ast}, \label{eq:global-convergence-growth-bound}
\endeq
the series \rref{eq:Fliess-operator-defined}
defines an operator from the extended space
$\Lpme (t_0)$ into $C[t_0, \infty)$,  where
\begin{align*}
\Lpme(t_0):=&
\{u:[t_0,\infty)\rightarrow \re^m:u_{[t_0,t_1]}\in \Lpm[t_0,t_1], \\
&\forall t_1 \in (t_0,\infty)\},
\end{align*}
and $u_{[t_0,t_1]}$ denotes the restriction of $u$ to $[t_0,t_1]$ \citep{Gray-Wang_SCL02}.
The set of all such {\em globally convergent} series is designated by
$\allseriesellGC$.

Given Fliess operators $F_c$ and $F_d$, where $c,d\in\allseriesell$,
the parallel and product connections satisfy $F_c+F_d=F_{c+d}$ and $F_cF_d=F_{c\shuffle d}$,
respectively \cite{Fliess_81}.
When Fliess operators $F_c$ and $F_d$ with $c\in\allseriesell$ and
$d\in\allseriesm$ are interconnected in a cascade fashion, the composite
system $F_c\circ F_d$ has the
Fliess operator representation $F_{c\circ d}$, where
the {\em composition product} of $c$ and $d$
is given by
\begdi
c\circ d=\sum_{\eta\in X^\ast} (c,\eta)\,\psi_d(\eta)(1)
\enddi
\cite{Ferfera_79,Ferfera_80}.
Here $\psi_d$ is the continuous (in the ultrametric sense) algebra homomorphism
from $\allseries$ to $\Endallseries$ uniquely specified by
$\psi_d(x_i\eta)=\psi_d(x_i)\circ \psi_d(\eta)$ with
\begdi \label{eq:psi-d-on-words}
\psi_d(x_i)(e)=x_0(d_i\shuffle e),
\enddi
$i=0,1,\ldots,m$
for any $e\in\allseries$,
and where $d_i$ is the $i$-th component series of $d$
($d_0:=1$).
$\psi_d(\emptyset)$ is the identity map on $\allseries$.
This composition product is associative and $\re$-linear in its left argument.
\begin{figure}[t]
\begin{center}
\includegraphics*[scale=0.45]{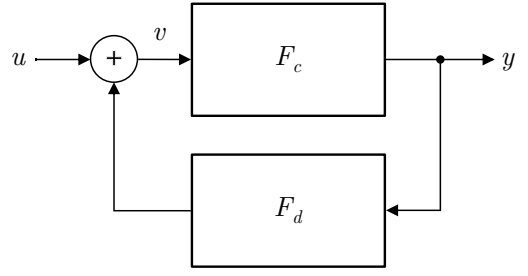}
\end{center}
\caption{Feedback connection}
\label{fig:feedback-with-v}
\end{figure}
In the event that two Fliess operators are interconnected
to form a feedback system as shown in Figure~\ref{fig:feedback-with-v},
it was shown in \cite{Gray-Li_05} that
there always exists a unique generating series $c@d$ such that
$y=F_{c@d}[u]$ whenever $c,d\in\allseriesmLC$.
This so called {\em feedback product} of $c$ and $d$
can be viewed
as the unique fixed point of a contractive
iterated map on a complete ultrametric space, but to compute it explicitly
requires Hopf algebraic tools such as those employed in
\cite{Gray-Duffaut_Espinosa_SCL11,Gray-Duffaut_Espinosa_FdB14} for
SISO systems. The multivariable case is considered in the next section.

\section{Hopf Algebra for Multivariable Output Feedback Group}
\label{sec:FdB-Hopf-algebra}

Consider the set of operators
$
\Fliessdelta:=\{I+F_c:c\in\allseriesm\},
$
where $I$ denotes the identity operator.
It is convenient to introduce the symbol
$\delta$ as the (fictitious) generating series for the identity map. That is,
$F_\delta:=I$ such that $I+F_c:=F_{\delta+c}=F_{c_\delta}$ with
$c_\delta:=\delta+c$.
The set of all such generating series for
$\Fliessdelta$ will be denoted by $\allseriesdeltam$.
The first theorem describes the multivariable output feedback group which is
at the heart of all the analysis in this paper. The group product
is described in terms of the
{\em modified composition product} of $c\in\allseriesell$ and $d\in\allseriesm$,
namely,
\begdi
c\modcomp d=\sum_{\eta\in X^\ast} (c,\eta)\, \phi_d(\eta)(1),
\enddi
where $\phi_d$ is the continuous (in the ultrametric sense) algebra homomorphism from
$\allseries$ to $\Endallseries$ uniquely specified by
$\phi_d(x_i\eta)=\phi_d(x_i)\circ \phi_d(\eta)$ with
\begdi \label{eq:phi-d-on-words}
\phi_d(x_i)(e)=x_ie+x_0(d_i\shuffle e),
\enddi
$i=0,1,\ldots,m$
for any $e\in\allseries$, and where $d_0:=0$. Again,
$\phi_d(\emptyset)$ is the identity map on $\allseries$ \cite{Gray-Li_05}.
It can be easily shown that $F_c\circ (I+F_d)=F_{c\modcomp d}$ and for any $x_i\in X$
\begeq
(x_ic)\modcomp d=x_i(c\modcomp d)+x_0(d_i\shuffle (c\modcomp d)). \label{eq:xic-cmod-d-identity}
\endeq
The following (non-associativity) identity was proved in \cite{Li_04}
\begeq
(c\modcomp d)\modcomp e=c\modcomp(d\modcomp e+e)\label{eq:cmod-non-associative}
\endeq
for all $c\in\allseriesell$ and $d,e\in\allseriesm$.
The lemma below will be also useful. Its proof is deferred until
Section~\ref{sec:FdB-Hopf-algebra}, when all the appropriate tools are
available.

\begle \label{le:c-modcomp-d-equals-K}
Let $d\in\allseriesm$ be fixed. Then $c\modcomp d=K\in\re^\ell$ if and only
if $c=K$.\footnote{For notational convenience, $c=K\emptyset$ is written as $c=K$.}
\endle

The central idea is that $(\Fliessdelta,\circ,I)$ forms a group of operators
under the composition
\begdi
F_{c_\delta}\circ F_{d_\delta}:=(I+F_c)\circ(I+F_d)
= F_{c_\delta\circ d_\delta},
\enddi
where $c_\delta\circ d_\delta:=\delta+d+c\modcomp d=:\delta+c\circledcirc d$.\footnote{The same symbol will be
used for composition on $\allseriesm$ and $\allseriesdeltam$.
As elements in these two sets have a distinct
notation, i.e., $c$ versus $c_\delta$, respectively, it will always be clear which product is at play.}$\,\,$%
Given the uniqueness of generating series of Fliess operators,
this assertion is equivalent to the following theorem.

\begin{table*}
\caption{Bases in the gradings $V=\bigoplus_{k\geq 0}V_k$ and $H=\bigoplus_{k\geq 0}H_k$. Here $i_\ell,j_\ell\neq 0$.}
\label{tbl:char-in-graded-H}
\begin{center}
\renewcommand{\arraystretch}{1.7}
\begin{tabular}{|c|c|c|c|c|} \hline
$k$ & $V_k$ & $H_k$ & $\dim(V_k)$ & $\dim(H_k)$ \\ \hline\hline
0 & $\mbf{1}$                   & $\mbf{1}$ & 1 & $1$ \\ \hline
1 & $a^{i_1}_\emptyset$       & $a^{i_1}_\emptyset$ & $m$ & $m$ \\[0.02in] \hline
2 & $a^{i_1}_{x_{j_1}}$           & $a^{i_1}_{x_{j_1}},a^{i_1}_{\emptyset}a^{i_2}_\emptyset$ & $m^2$ & $2m^2$ \\[0.02in] \hline
3 & $a^{i_1}_{x_0},a^{i_1}_{x_{j_1}x_{j_2}}$ & $a^{i_1}_{x_0},a^{i_1}_{x_{j_1}x_{j_2}},a^{i_1}_\emptyset a^{i_2}_\emptyset a^{i_3}_\emptyset,a^{i_1}_\emptyset a^{i_2}_{x_{j_1}}$
& $m+m^3$ & $m+3m^3$ \\[0.02in] \hline
4 & $a^{i_1}_{x_0x_{j_1}},a^{i_1}_{x_{j_1}x_0},a^{i_1}_{x_{j_1}x_{j_2}x_{j_3}}$ &
\begin{tabular}{c}
$a^{i_1}_{x_0x_{j_1}},a^{i_1}_{x_{j_1}x_0},a^{i_1}_{x_{j_1}x_{j_2}x_{j_3}},a^{i_1}_\emptyset a^{i_2}_\emptyset a^{i_3}_\emptyset a^{i_4}_\emptyset,$ \\
$a^{i_1}_\emptyset a^{i_2}_\emptyset a^{i_3}_{x_{j_1}},a^{i_1}_{\emptyset}a^{i_2}_{x_0},a^{i_1}_{\emptyset}a^{i_2}_{x_{j_1}x_{j_2}},a^{i_1}_{x_{j_1}}a^{i_2}_{x_{j_2}}$
\end{tabular}
& $2m^2+m^4$ & $3m^2+5m^4$ \\[0.15in] \hline
\end{tabular}
\end{center}
\end{table*}

\begth \label{th:allseriesdeltam-is-group}
The triple $(\allseriesdeltam,\circ,\delta)$ is a group.
\endth

\begpr
By design, $\delta$ is the identity element of the group. The associativity of the
product can be established in a manner similar to the SISO case addressed
in \cite{Gray-Duffaut_Espinosa_SCL11}. (See \cite{Gray_MTNS14} for an
alternative approach.) The existence of an inverse will
be handled differently here (more directly) via Lemma~\ref{le:c-modcomp-d-equals-K}.
Specifically, for a fixed $c_\delta\in\allseriesdeltam$, the composition inverse, $c_\delta^{-1}=\delta+c^{-1}$, must satisfy
$c_\delta\circ c_\delta^{-1}=\delta$ and $c_\delta^{-1}\circ c_\delta=\delta$,
which reduce, respectively, to
\begin{subequations}
\begin{align}
c^{-1}&=(-c)\modcomp c^{-1} \label{eq:cdelta-right-inverse} \\
c&=(-c^{-1})\modcomp c. \label{eq:cdelta-left-inverse}
\end{align}
\end{subequations}
It was shown in \cite{Gray-Li_05}
that $e\mapsto (-c)\modcomp e$
is always a contraction in the ultrametric sense on $\allseriesm$ as
a complete ultrametric space and thus has a unique fixed point.
So it follows directly that $c_\delta^{-1}$ is a right inverse of $c_\delta$, i.e., satisfies
\rref{eq:cdelta-right-inverse}.
\markblue{To see that this same series is also a left inverse, first
observe that \rref{eq:cdelta-right-inverse} is equivalent to
\begeq \label{eq:left-inverse-RHS-version}
c^{-1}\modcomp 0 +c\modcomp c^{-1}=0,
\endeq
using the identity $c^{-1}\modcomp 0=c^{-1}$ and the left linearity of the
modified composition product. Substituting \rref{eq:left-inverse-RHS-version} back into
itself where zero appears and applying \rref{eq:cmod-non-associative}
gives
\begin{align*}
c^{-1}\modcomp (c\modcomp c^{-1}+c^{-1})+c\modcomp c^{-1}&=0 \\
(c^{-1}\modcomp c)\modcomp c^{-1}+c\modcomp c^{-1}&=0.
\end{align*}
Again from left linearity of the modified composition product it follows that
\begdi
(c^{-1}\modcomp c+c)\modcomp c^{-1}=0.
\enddi
Finally, Lemma~\ref{le:c-modcomp-d-equals-K} implies that $c^{-1}\modcomp c+c=0$,}
which is equivalent to \rref{eq:cdelta-left-inverse}. This concludes the proof.
\endpr

A Fa\`{a} di Bruno type Hopf algebra is now defined for the
output feedback group.
The coordinate maps for this algebra have the form
\begdi \label{eq:character-maps}
a^i_\eta:\allseriesm\rightarrow \re:c\mapsto (c_i,\eta),
\enddi
where $\eta\in X^\ast$ and $i=1,2,\ldots,m$.\footnote{The use of this terminology will be justified in the proof of Lemma~\ref{le:antipode-is-group-inverse}.}
Let $V$ denote the $\re$-vector space spanned by these maps.
If the {\em degree} of $a^i_{\eta}$ is defined as $\deg(a^i_{\eta})=2\abs{\eta}_{x_0}+\sum_{j=1}^m\abs{\eta}_{x_j}+1$,
then $V$ is a connected graded vector space, that is, $V=\bigoplus_{n\geq 0} V_n$
with
\begdi
V_n=\spanset_\re\{a_\eta^i:\deg(a^i_\eta)=n\}, \;\;n>0,
\enddi
$V_0=\re\mbf{1}$, and $\mbf{1}$ maps every $c\in\allseriesm$ to 1.

Consider next the free unital commutative $\re$-algebra, $H$, with product
\begdi \label{eq:mu-product}
\mu:a^i_\eta\otimes a^j_{\xi}\mapsto a^i_{\eta}a^j_{\xi}
\enddi
and \markblue{unit $\mbf{1}$.}
This product is clearly associative.
The graduation on $V$ induces a connected graduation on $H$
with
$\deg(a^i_\eta a^j_\xi)=\deg(a^i_\eta)+\deg(a^j_\xi)$ and $\deg(\mbf{1})=0$.
Specifically,
$H=\bigoplus_{n\geq 0} H_n$, where $H_n$ is the set of
all elements of degree $n$ and $H_0=\re\mbf{1}$.
Bases for these subspaces are given in Table~\ref{tbl:char-in-graded-H}.

Three coproducts are now introduced.
\markblue{The first coproduct is used to define the Hopf algebra on $H$.
The remaining two coproducts provide a recursive manner in which to compute it.
Recalling that $c_\delta\circ d_\delta=\delta+c\circledcirc d$,
define $\Delta$ for any $a^i_\eta\in V_+$ such that
\begdi
\Delta a^i_\eta(c,d)=a^i_\eta(c\circledcirc d)=(c_i\circledcirc d,\eta).
\enddi}%
The coassociativity of $\Delta$ follows from the associativity of the product
$c\circledcirc d$ \cite{Gray-Duffaut_Espinosa_SCL11}.
Specifically, for any $c,d,e\in\allseriesm$:
\begin{align*}
(\id\otimes\Delta)\circ \Delta a^i_\eta(c,d,e)&=(c_i\circledcirc(d\circledcirc e),\eta)\\
&=((c\circledcirc d)_i\circledcirc e,\eta) \\
&=(\Delta\otimes \id)\circ \Delta a^i_\eta(c,d,e).
\end{align*}
Therefore, $(\id\otimes\Delta)\circ \Delta=(\Delta\otimes \id)\circ \Delta$ as required.
\markblue{The following result motivates the primary interest in proving that
$(H,\mu,\Delta)$ is a Hopf algebra, namely, that it has
an antipode corresponding to the group
inverse for $(\allseriesdeltam,\circ,\delta)$.}

\markblue{
\begle \label{le:antipode-is-group-inverse}
The Hopf algebra $(H,\mu,\Delta)$ has an antipode $S$ satisfying
$a_\eta^i(\gsc^{-1})=(Sa^i_\eta)(\gsc)$ for all $\eta\in X^\ast$ and $c\in\allseriesm$.
\endle
\begpr
First observe that for
each $\gsc_\delta\in\allseriesdeltam$, one can identify a character map $\Phi_\gsc\in L(H,\re)$ as
\begdi
	\Phi_\gsc:a^i_\eta\mapsto a^i_\eta(\gsc)=(\gsc_i,\eta),
\enddi
where by
$\Phi_\gsc(\mbf{1})=1$ and
\begdi
	\Phi_\gsc(a^i_\eta a^j_\xi)=a^i_\eta(\gsc) a^j_\xi(\gsc)=\Phi_\gsc(a^i_\eta)\Phi_\gsc(a^j_\xi).
\enddi
Coordinate maps, therefore, should satisfy \rref{eq:general-coordinate-maps}, specifically,
\begin{align*}
	(\Phi_\gsc \star \Phi_\gsd)(a^i_\eta)&=\mu\circ(\Phi_\gsc\otimes \Phi_\gsd)\circ\Delta a^i_\eta \\
		&=\sum \Phi_\gsc\left({a^i_\eta}_{(1)}\right)\Phi_\gsd\left({a^i_\eta}_{(2)}\right) \\
		&=\sum {a^i_\eta}_{(1)}(\gsc) {a^i_\eta}_{(2)}(\gsd) \\
		&= \Delta a^i_\eta(\gsc,\gsd) \\
		&= a_\eta^i(\gsc\circledcirc \gsd) \\
        &= a_{\eta}^i(c_\delta\circ d_\delta),
\end{align*}
where the summation is taken over all terms that appear in
$\Delta a_\eta^i$
(following the notation of Sweedler \cite{Sweedler_69}).
From this identification between the convolution of characters and the group product on $\allseriesdeltam$,
it is clear that $\Phi^{\star -1}_c=\Phi_{c^{-1}}$.
Using \rref{eq:convinv-via-antipode} then
$
\Phi_{c^{-1}}=\Phi_\gsc^{\star -1}=\Phi_\gsc \circ S,
$
where $\gsc_\delta^{-1}=\delta+\gsc^{-1}$, so that
$
a_\eta^i(\gsc^{-1})=(Sa^i_\eta)(\gsc)
$
as desired.
\endpr
}

The second coproduct is
$\Delta_\shuffle^j(V_+)\subset V_+\otimes V_+$, which is isomorphic to $\sh^\ast$ via the
coordinate maps.
That is,
\begin{subequations}
\label{eq:shuffle-coproduct-induction}
\begin{align}
\Delta_{\shuffle}^ja^i_{\emptyset}&=a^i_{\emptyset}\otimes a^j_{\emptyset} \\
\Delta_{\shuffle}^j\circ\theta_k&=(\theta_k\otimes \id+\id\otimes \theta_k)\circ\Delta_{\shuffle}^j,
\end{align}
\end{subequations}
where
$\id$ is the identity map on $V_+$, and
$\theta_k$ denotes the endomorphism on $V_+$ specified by
$\theta_ka^i_\eta=a^i_{x_k\eta}$ for $k=0,1,\ldots,m$ and $i,j=1,2,\ldots,m$.

\begex \label{ex:character-shuffle-coproduct}
The first few terms of $\Delta_{\shuffle}^j$ are:
\begin{align*}
\Delta_{\shuffle}^ja^i_{\emptyset}&=a^i_{\emptyset}\otimes a^j_{\emptyset} \\
\Delta_{\shuffle}^ja^i_{x_{i_1}}&=a^i_{x_{i_1}}\otimes a^j_{\emptyset}+a^i_{\emptyset}\otimes a^j_{x_{i_1}} \\
\Delta_{\shuffle}^ja^i_{x_{i_2}x_{i_1}}&=a^i_{x_{i_2}x_{i_1}}\otimes a^j_{\emptyset}+a^i_{x_{i_2}}\otimes a^j_{x_{i_1}}+a^i_{x_{i_1}}\otimes a^j_{x_{i_2}}+ \\
&\hspace*{0.18in}a^i_{\emptyset}\otimes a^j_{x_{i_2}x_{i_1}} \\
\Delta_{\shuffle}^ja^i_{x_{i_3}x_{i_2}x_{i_1}}&=a^i_{x_{i_3}x_{i_2}x_{i_1}}\otimes a^j_{\emptyset}+a^i_{x_{i_3}x_{i_2}}\otimes a^j_{x_{i_1}}+ \\
&\hspace*{0.18in}
a^i_{x_{i_3}x_{i_1}}\otimes a^j_{x_{i_2}}+a^i_{x_{i_3}}\otimes a^j_{x_{i_2}x_{i_1}}+\\
&\hspace*{0.18in}a^i_{x_{i_2}x_{i_1}}\otimes a^j_{x_{i_3}}+a^i_{x_{i_2}}\otimes a^j_{x_{i_3}x_{i_1}}+ \\
&\hspace*{0.18in}
a^i_{x_{i_1}}\otimes a^j_{x_{i_3}x_{i_2}}+a^i_{\emptyset}\otimes a^j_{x_{i_3}x_{i_2}x_{i_1}}.
\end{align*}
\endex

The third coproduct is \markblue{$\tilde{\Delta}a_\eta^i=\Delta a_\eta^i-\mbf{1}\otimes a_\eta^i$
or, equivalently,
the coproduct induced by the identity}
\begdi
\tilde{\Delta}a^i_{\eta}(c,d)=(c_i\modcomp d,\eta)
=\sum a^i_{\eta(1)}(c) a^i_{\eta(2)}(d). \label{eq:tilde-delta-identity}
\enddi
A key observation is that this coproduct
can be computed recursively as described in
the next lemma, which is the multivariable version of
Proposition 3 in \cite{Foissy_13}.
\markblue{It is not difficult to show using (2) and (3) of this lemma
that $a^i_{\eta(1)}\in V_+$ and $a^i_{\eta(2)}\in H$, and thus, $\tilde{\Delta}V_+\subseteq V_+\otimes H$.}

\begle \label{le:tilde-delta-inductions}
The following identities hold:
\begdes
\item[\hspace*{0.15in}(1)] $\tilde{\Delta}a^i_\emptyset=a^i_{\emptyset}\otimes \mbf{1}$
\item[\hspace*{0.15in}(2)]
$\tilde{\Delta}\circ \theta_i=
(\theta_i\otimes \id)\circ \tilde{\Delta}$
\item[\hspace*{0.15in}(3)]
$\tilde{\Delta}\circ \theta_0=(\theta_0\otimes \id)\circ \tilde{\Delta}+
(\theta_i\otimes \mu)\circ(\tilde{\Delta}\otimes \id)\circ \Delta_{\shuffle}^i$,
\enddes
$i=1,2,\ldots,m$,
where $\id$ denotes the identity map on $H$.\footnote{The Einstein summation notation
is used in item (3) and throughout
to indicate summations from either 0 or 1 to $m$, e.g., $\sum_{i=1}^m a_ib^i=a_ib^i$.
It will be clear from the context
which lower bound is applicable.}
\endle

\begpr $\,$\\
\noindent(1) First note that any series $c$ can be uniquely decomposed as
$c=(c,\emptyset)\emptyset+ x_ic^i$, $i=0,1,\ldots,m$,
where the series $c^i$ are arbitrary.
In which case, using the left linearity of the modified composition product
and \rref{eq:xic-cmod-d-identity}, it follows that
\begin{align*}
\tilde{\Delta}a^i_\emptyset(c,d)&=a^i_{\emptyset}(c\modcomp d)
=a^i_{\emptyset}\left((c,\emptyset)\emptyset+ (x_jc^j)\modcomp d \right) \\
&=(c_i,\emptyset)+a^i_\emptyset(x_j(c^j\modcomp d)+x_0(d_j\shuffle (c^j\modcomp d))) \\
&=(c_i,\emptyset)=(a^i_\emptyset\otimes \mbf{1})(c,d).
\end{align*}

\noindent (2) For any $\eta\in X^\ast$ observe
\begin{align*}
(\tilde{\Delta} \circ \theta_i) a^j_\eta (c,d)&= \tilde{\Delta} a^j_{x_i\eta}(c,d) \\
&= a^j_{x_i\eta}(x_k(c^k\modcomp d)+x_0(d_k\shuffle (c^k\modcomp d))) \\
&= a^j_\eta(c^i\modcomp d) \\
&= \tilde{\Delta}a^j_\eta(c^i,d) \\
&=\sum a^j_{\eta(1)}\otimes a^j_{\eta(2)}(c^i,d) \\
&=\sum \theta_i(a^j_{\eta(1)})\otimes a^j_{\eta(2)}(c,d) \\
&=(\theta_i\otimes \id)\circ \tilde{\Delta}a^j_{\eta} (c,d).
\end{align*}
Note that since $a^j_{\eta(1)}\in V_+$, the operation $\theta_i(a^j_{\eta(1)})$ is well defined.

\noindent (3) Proceeding as in the previous item, it follows that
\begin{align*}
\lefteqn{(\tilde{\Delta} \circ \theta_0) a^i_\eta(c,d)}\\
&=a^i_{x_0\eta}(c\modcomp d) \\
&= a^i_{x_0\eta}(x_j(c^j\modcomp d)+x_0(d_j\shuffle (c^j\modcomp d))) \\
&= a^i_{\eta}(c^0\modcomp d+d_j\shuffle (c^j\modcomp d)) \\
&= a^i_{\eta}(c^0\modcomp d)+\sum_{j=1}^m \Delta_{\shuffle}^j a^i_\eta(c^j\modcomp d,d) \\
&= a^i_{\eta}(c^0\modcomp d)+
\sum_{j=1}^m \sum_{\xi,\nu\in X^\ast} (\eta,\xi\shuffle \nu)\,
a^i_{\xi}(c^j\modcomp d)a^j_\nu(d) \\
&= \tilde{\Delta}a^i_{\eta}(c^0,d)+ \sum_{j=1}^m \sum_{\xi,\nu\in X^\ast} (\eta,\xi\shuffle \nu)\,
(\tilde{\Delta}a^i_{\xi}\otimes a^j_\nu)(c^j,d,d) \\
&= (\theta_0\otimes \id)\circ \tilde{\Delta}a^i_{\eta}(c,d)+(\theta_j\otimes \id)\circ \\
&\hspace*{0.18in} \sum_{\xi,\nu\in X^\ast} (\eta,\xi\shuffle \nu)\,
(\tilde{\Delta}a^i_{\xi}\otimes a^j_\nu)(c,d,d) \\
&= (\theta_0\otimes \id)\circ \tilde{\Delta}a^i_{\eta}(c,d)+ (\theta_j\otimes \mu)
\circ(\tilde{\Delta}\circ \id)\circ \Delta_{\shuffle}^ja^i_\eta (c,d).
\end{align*}
\endpr

The next theorem is a central result of the paper.

\begth
$(H,\mu,\Delta)$ is a connected graded commutative unital Hopf algebra.
\endth

\begpr
From the development above, it is clear that $(H,\mu,\Delta)$
is a bialgebra \markblue{with unit $\mbf{1}$ and counit $\varepsilon$ defined by $\varepsilon(a_{\eta})=0$ for all
$\eta\in X^\ast$ and $\varepsilon (\mbf{1})=1$ (see also \cite[equation (14)]{Gray-Duffaut_Espinosa_SCL11}).}
Here it is shown that this bialgebra is graded and connected.
Therefore, $H$ automatically has an antipode, and thus, is a Hopf algebra \cite{Figueroa-Gracia-Bondia_05}.
Specifically, since the algebra $H$ is graded by $H_n$, $n\geq 0$ with $H_0=\re\mbf{1}$,
it only needs to be shown for any $a^i_\eta\in V_+$
that
\begeq
\tilde{\Delta}a^i_\eta\in (V_+\otimes H)_{n}:=\bigoplus_{j+k=n\atop j\geq 1,k\geq 0} V_j\otimes  H_k. \label{eq:tilde-delta-grading}
\endeq
This fact is evident from the first few terms computed via Lemma~\ref{le:tilde-delta-inductions}:
\begin{align*}
n=1&:\tilde{\Delta} a_{\emptyset}^i= a_{\emptyset}^i\otimes \mbf{1}\\
n=2&:\tilde{\Delta} a_{x_j}^i= a_{x_j}^i\otimes \mbf{1} \\
n=3&:\tilde{\Delta} a_{x_0}^i=a_{x_0}^i\otimes \mbf{1}+ a_{x_\ell}^i \otimes a_{\emptyset}^\ell \\
n=3&:\tilde{\Delta} a_{x_jx_k}^i=a_{x_jx_k}^i\otimes \mbf{1} \\
n=4&:\tilde{\Delta} a_{x_0x_j}^i= a_{x_0x_j}^i\otimes \mbf{1}+ a_{x_\ell}^i \otimes a_{x_j}^\ell + a_{x_\ell x_j}^i \otimes a_{\emptyset}^\ell \\
n=4&:\tilde{\Delta} a_{x_jx_0}^i=  a_{x_jx_0}^i\otimes \mbf{1}+a_{x_jx_\ell}^i \otimes a_{\emptyset}^\ell \\
n=4&:\tilde{\Delta} a_{x_jx_kx_l}^i=a_{x_jx_kx_l}^i\otimes \mbf{1} \\
n=5&:\tilde{\Delta} a_{x_0^2}^i= a_{x_0^2}^i\otimes \mbf{1}+ a_{x_\ell}^i \otimes a_{x_0}^\ell + a_{x_\ell x_0}^i \otimes a_{\emptyset}^\ell +  \\
&\hspace*{0.65in}  a_{x_0x_\ell}^i \otimes a_{\emptyset}^\ell +a_{x_\ell x_\nu}^i \otimes a_{\emptyset}^\ell a_{\emptyset}^\nu,
\end{align*}
where $i,j,k,l=1,2,\ldots m$.
In which case, using the identities $\Delta(a^i_\eta a^j_\xi)=\Delta a_\eta^i \Delta a_\xi^j$
and $\Delta a_\eta^i=\tilde{\Delta}a_\eta^i+\mbf{1}\otimes a_\eta^i$, it follows that
$
\Delta H_n\subseteq (H\otimes H)_n,
$
and this would complete the proof. To prove \rref{eq:tilde-delta-grading},
the following facts are essential:
\begen
\item $\deg(\theta_la^i_\eta)=\deg(a^i_\eta)+1$, $l=1,2,\ldots,m$
\item $\deg(\theta_0a^i_\eta)=\deg(a^i_\eta)+2$
\item $\Delta_\shuffle^j a^i_\eta\in (V_+\otimes V_+)_{n+1}$, $n=\deg(a_\eta^i)$.
\enden
The proof is via
induction on the length of $\eta$. When $\abs{\eta}=0$ then clearly
$\tilde{\Delta}a^i_\emptyset=a^i_\emptyset\otimes \mbf{1}\in V_1\otimes H_0$ and $n=1$.
Assume now that \rref{eq:tilde-delta-grading} holds for words up to some fixed length
$\abs{\eta}\geq 0$. Let $n=\deg(a_\eta^i)$.
There are two ways to increase the
length of $\eta$.
First consider  $a^i_{x_l\eta}$ for some $l\neq 0$.
From item 1 above $\deg(a^i_{x_l\eta})=n+1$, and from Lemma~\ref{le:tilde-delta-inductions}
$\tilde{\Delta}a^i_{x_l\eta}=(\theta_l\otimes \id)\circ \tilde{\Delta}a^i_\eta$.
Therefore, using the induction hypothesis, $\tilde{\Delta}a^i_{x_la_\eta}\in \bigoplus_{j+k=n}
V_{j+1}\otimes H_k\subset (V\otimes H)_{n+1}$, which proves the assertion.
Consider next $a^i_{x_0\eta}$. From item~2 above $\deg(a^i_{x_0\eta})=n+2$.
Lemma~\ref{le:tilde-delta-inductions} is employed as in the first case.
First note that from item 3 above
$\Delta_\shuffle^j a^i_\eta\in (V_+\otimes V_+)_{n+1}$, and so
using the induction hypothesis it follows that
$(\tilde{\Delta}\otimes \id)\circ \Delta_\shuffle^j a^i_\eta\in (V_+\otimes H\otimes V_+)_{n+1}$. In which case,
$(\theta_i\otimes \mu)\circ (\tilde{\Delta}\otimes \id)\circ \Delta_\shuffle^j a^i_\eta\in (V_+\otimes H)_{n+2}$. By a similar argument,
$(\theta_0\otimes \id)\circ \tilde{\Delta} a_\eta^i\in (V_+\otimes H)_{n+2}$.
Thus,
$\tilde{\Delta}a^i_{x_0\eta}\in (V_+\otimes H)_{n+2}$, which again proves the assertion
and completes the proof.
\endpr

The deferred proof from Section~\ref{sec:preliminaries} is presented next.

\noindent{\sc Proof of Lemma~\ref{le:c-modcomp-d-equals-K}.}
The only non trivial claim is that $c\modcomp d=K$ implies $c=K$.
If $c\modcomp d=K$ then clearly
$K_i=a^i_{\emptyset}(c\modcomp d)=\tilde{\Delta}a^i_{\emptyset}(c,d)=a^i_{\emptyset}c$,
$i=1,2,\ldots,\ell$.
Furthermore, for any $x_j\in X$ with $j\neq 0$, $0=a^i_{x_j}(c\modcomp d)=\tilde{\Delta}a^i_{x_j}(c,d)=a^i_{x_j}c$,
$i=1,2,\ldots,\ell$.
Now suppose $a^i_\eta c=0$, $i=1,2,\ldots,\ell$ for all $a^i_{\eta}\in V_k$ with $k=1,2,\ldots,n$.
Then for any $x_j\in X$
\begdi
0=\tilde{\Delta}a^i_{x_j\eta}(c,d)=a^i_{x_j\eta}c+
\sum_{a^i_{x_j\eta(2)}\neq 1} a^i_{x_j\eta(1)}(c)\;a^i_{x_j\eta(2)}(d),
\enddi
where in general $a^i_{x_j\eta(1)}\neq a^i_{\emptyset}$. Therefore,
$a^i_{x_j\eta}c=0$, $i=1,2,\ldots,\ell$. In which, case $c=K$.

The antipode of any graded connected Hopf algebra can be computed as described in
the following theorem. It can be viewed as being {\em partially} recursive in that
the coproduct needs to be computed first before the antipode recursion can be applied.

\begth {\rm \cite{Figueroa-Gracia-Bondia_05}} \label{th:antipode-induction}
The antipode, $S$, of any graded connected Hopf algebra $(H,\mu,\Delta)$
can be computed for any $a\in H_{k}$, $k\geq 1$ by
\begdi
S a=-a-\sum (S a^\prime_{(1)})a^\prime_{(2)}=-a-\sum  a^\prime_{(1)}S a^\prime_{(2)},
\label{eq:general-antipode-induction}
\enddi
where the reduced coproduct is
$\Delta^\prime a=\Delta a-a\otimes \mbf{1}-\mbf{1}\otimes a=\sum a^\prime_{(1)}a^\prime_{(2)}$.
\endth

The next theorem provides a {\em fully} recursive algorithm to compute the antipode for the
output feedback group.

\begth
The antipode, $S$, of any $a_\eta^i\in V_+$ for the output feedback group can be computed by the following algorithm:
\begdes
\item[\hspace*{0.24in}i.] Recursively compute $\Delta_{\shuffle}^j$ via \rref{eq:shuffle-coproduct-induction}.
\item[\hspace*{0.2in}ii.] Recursively compute $\tilde{\Delta}$ via Lemma~\ref{le:tilde-delta-inductions}.
\item[\hspace*{0.17in}iii.] Recursively compute $S$ via Theorem~\ref{th:antipode-induction} with
$\Delta^\prime a_\eta^i=$ \\ \hspace*{0.2in} $\tilde{\Delta}a_\eta^i-a_\eta^i\otimes \mbf{1}$.
\enddes
\endth
\begpr
In light of the previous results, the only detail is the minor observation
that $S$ is the antipode of
the Hopf algebra with coproduct $\Delta a=\tilde{\Delta}a+\mbf{1}\otimes a$.
In which case, the corresponding
reduced coproduct is as described in step iii.
\endpr

Applying the algorithm above gives the antipode of the first few
coordinate maps:
\begin{align*}
H_1&:S a^i_\emptyset=-a^i_\emptyset \\
H_2&:S a^i_{x_j}=-a^i_{x_j} \\
H_3&:S a^i_{x_0}=-a^i_{x_0}+ a^i_{x_\ell}a^\ell_\emptyset \\
H_3&:S a^i_{x_jx_k}=-a^i_{x_jx_k} \\
H_4&:S a^i_{x_0x_j}=-a^i_{x_0x_j}+a^i_{x_\ell}a^\ell_{x_j}+a^i_{x_\ell x_j}a^\ell_\emptyset \\
H_4&:S a^i_{x_jx_0}=-a^i_{x_jx_0}+ a^i_{x_jx_\ell}a^\ell_{\emptyset} \\
H_4&:S a^i_{x_j x_k x_l}=-a^i_{x_j x_k x_l} \\
H_5&:S a^i_{x_0^2}=-a^i_{x_0^2}-(S a^i_{x_\ell})a^\ell_{x_0}-(S a^i_{x_\ell x_0})a^\ell_\emptyset - \\
&\hspace*{0.58in}(S a^i_{x_0x_\ell})a^\ell_\emptyset - (S a^i_{x_\ell x_\nu})a^\ell_\emptyset a^\nu_\emptyset \\
&\hspace*{0.4in}=-a^i_{x_0^2}-(-a^i_{x_\ell})a^\ell_{x_0}-(-a^i_{x_\ell x_0}+\cancel{a^i_{x_\ell x_\nu}a^\nu_\emptyset})a^\ell_\emptyset - \\
&\hspace*{0.58in}(-a^i_{x_0 x_\ell}+a^i_{x_\nu}a^\nu_{x_\ell}+a^i_{x_\nu x_\ell}a^\nu_{\emptyset})a^\ell_\emptyset- \\
&\hspace*{0.58in}\cancel{(-a^i_{x_\ell x_\nu})a^\ell_\emptyset a^\nu_\emptyset} \\
&\hspace*{0.4in}=-a^i_{x_0^2}+a^i_{x_\ell}a^\ell_{x_0}+a^i_{x_\ell x_0}a^\ell_\emptyset+ a^i_{x_0x_\ell}a^\ell_\emptyset- \\
&\hspace*{0.6in}a^i_{x_\nu}a^{\nu}_{x_\ell} a^\ell_{\emptyset}-a^i_{x_\nu x_\ell}a^\nu_{\emptyset}a^\ell_{\emptyset},
\end{align*}
where $i,j,k,l=1,2,\ldots m$.
The explicit calculations for $S a^i_{x_0^2}$ are shown above to demonstrate that this approach, not unexpectedly, involves some inter-term
cancellation. This is consistent with what is known about the classical Fa\`{a} di Bruno Hopf algebra and the Zimmermann
formula, which provides a cancellation free approach to computing the antipode
\cite{Einziger_10,Haiman-Schmitt_89}. This hints at the possibility of even more
efficient antipode algorithms, but this topic will not be pursued here.
It should also be noted that when $m=1$, i.e., the
SISO case, all the summations above vanish, and the identities reduce to those given
in \cite{Gray-Duffaut_Espinosa_SCL11}.

\begex
\label{ex:c-inverse-finite-Lie-rank}
{\rm
Suppose $c\in\allseriesmLC$ has finite {\em Lie rank},
that is, the range space of the Hankel mapping for $c$ defined on
the $\re$-vector space of Lie polynomials is finite. Then
$I+F_c$ has a finite dimensional control-affine state space
realization of the form
\begin{align*}
\dot{z}&= g_0(z)+\sum_{j=1}^m g_j(z)u_j,\;\;z(0)=z_0 \\
y_i&= h_i(z)+u_i,\;\;i=1,2,\ldots,m,
\end{align*}
where each $g_j$ and $h_i$ is an analytic vector field and function, respectively,
on some neighborhood $W$ of $z_0$ \cite{Fliess_83,Isidori_95}.
In which case,
\begeq
(c_i,\eta)=L_{g_{\eta}}h_i(z_0),\;\;\forall\eta\in X^\ast, \label{eq:c-equals-Lgh}
\endeq
where
\begdi
L_{g_{\eta}}h_i:=L_{g_{j_1}}\cdots L_{g_{j_k}}h_i,
\;\;\eta=x_{j_k}\cdots x_{j_1},
\enddi
the {\em Lie derivative} of $h_i$ with respect to $g_j$ is defined as
\begdi
L_{g_j} h_i: W\rightarrow \re: z \mapsto \frac{\partial h_i}{\partial z}(z)\, g_j(z),
\enddi
and $L_{g_\emptyset}h_i=h_i$.
\begin{figure}[t]
\begin{center}
\includegraphics[scale=0.45]{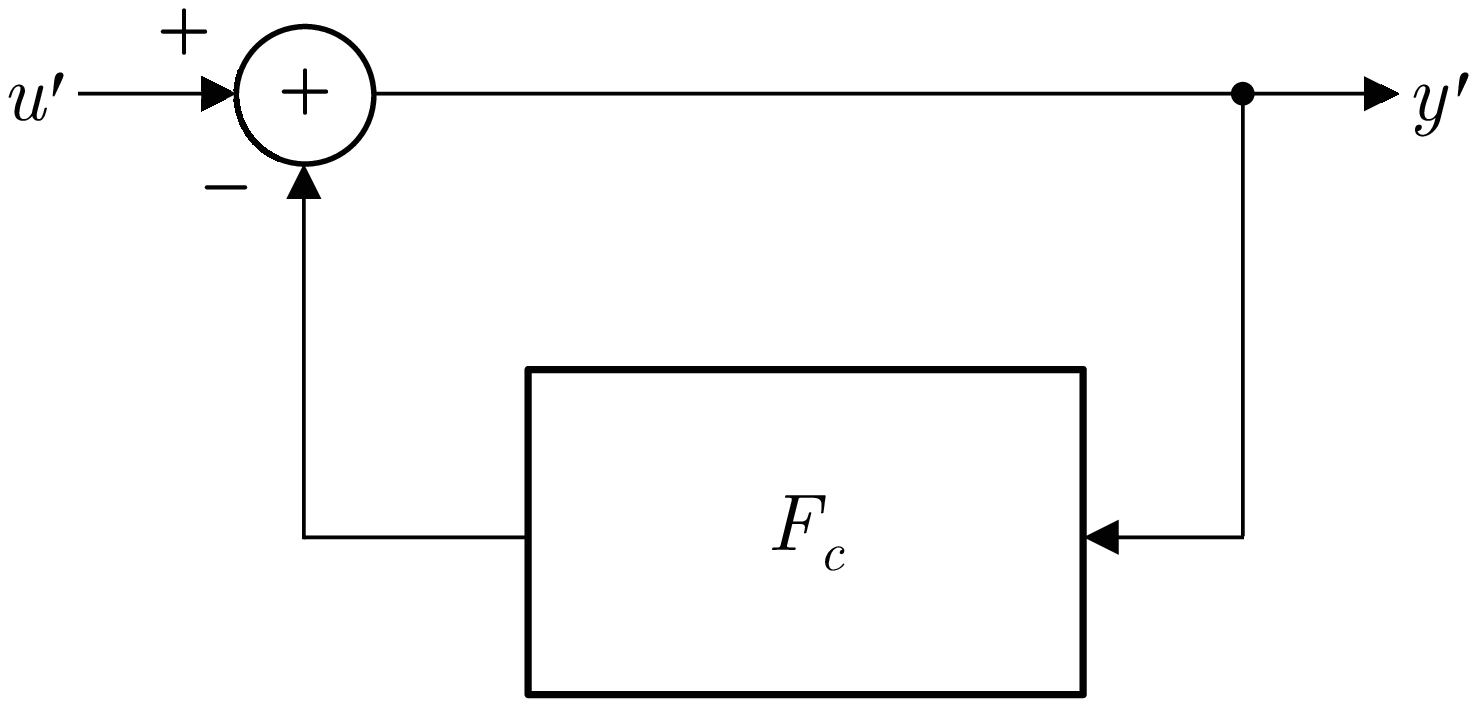}
\end{center}
\caption{Feedback implementation of the composition inverse of $I+F_c$.}
\label{fig:composition-inverse-IplusFc}
\end{figure}
It is not difficult to see that
the composition inverse of the return difference operator $I+F_c$, that is,
$(I+F_c)^{-1}=I+F_{c^{-1}}:u^\prime\mapsto y^\prime$,
is described by
the feedback system in Figure~\ref{fig:composition-inverse-IplusFc}.
A straightforward calculation gives a realization for $F_{c_i^{-1}}$,
namely, $(\{g_0-\sum_{j=1}^m g_jh_j,g_1,\ldots,g_m\},-h_i,z_0)$. Using this realization and \rref{eq:c-equals-Lgh},
it can be readily verified that Lemma~\ref{le:antipode-is-group-inverse} holds. For example,
\begin{align*}
(c_i^{-1},x_0)&=L_{g_0-\sum_j g_jh_j}(-h_i)(z_0) \\
&= -L_{g_0}h_i(z_0)+\sum_{j=1}^m (L_{g_j}h_i(z_0))h_j(z_0) \\
&= -(c_i,x_0)+\sum_{j=1}^m (c_i,x_j)(c_j,\emptyset) \\
&= (-a^i_{x_0} + a^i_{x_j}a^j_{\emptyset})c \\
&= (S a^i_{x_0})c.
\end{align*}
In the special case of a linear time-invariant system
with strictly proper $m\times m$ transfer function $H(s)$ and
state space realization $(A,B,C)$, the
corresponding components of the linear
generating series are $c_i=\sum_{k\geq 0} \sum_{j=1}^m(c_i,x_0^kx_j)\,x_0^kx_j$,
where $(c_i,x_0^kx_j)=C_iA^kB_j$, $k\geq 0$, and $C_i$, $B_j$ denote the $i$-th row of $C$ and the $j$-th column
of $B$, respectively. The composition inverse of the return difference matrix $I+H(s)$ is computed directly
as
\begdi
(I+C(sI-A)^{-1}B)^{-1}=I-C(sI-(A-BC))^{-1}B.
\enddi
Therefore, it
follows that
\begdi
(c_i^{-1},x_0^kx_j)=-C_i(A-BC)^kB_j,\;\;k\geq 0,\;\;i,j=1,2,\ldots,m.
\enddi
Expanding this product gives the expected antipode formulas. For example,
\begin{align*}
(c_i^{-1},x_0x_j)&= -C_i(A-BC)B_j \\
&= -C_iAB_j +C_iBCB_j \\
&=-C_iAB_j+\sum_{\ell=1}^m C_iB_\ell C_\ell B_j \\
&=-(c_i,x_0x_j)+\sum_{\ell=1}^m (c_i,x_\ell)(c_\ell,x_j) \\
&=(-a^i_{x_0x_j}+a^i_{x_\ell}a^\ell_{x_j}+a^i_{x_\ell x_j}a^\ell_\emptyset)c \\
&=(S a^i_{x_0x_j})c,
\end{align*}
where the fact that $(c,x_\ell x_j)=(c,\emptyset)=0$ has been used in the second to the last line.
It is worth repeating that the antipode formulas derived at the beginning of this
section required no state space setting. Hence, they still apply even when $c$ does {\em not}
have finite Lie rank.
}
\endex

The next theorem establishes that local convergence is preserved by
the composition inverse operation. This fact was proved for the SISO case
in \cite{Gray-Duffaut_Espinosa_SCL11} using only a grading of $H$. But here a
different approach is taken, one that
produces the exact radius of convergence for the
operation.

\begth \label{th:ROC-for-cinv_LC}
For any $c\in\allseriesmLC$ with
growth constants $K_c,M_c>0$ it follows that
\begeq
|(c^{-1},\eta)|\leq K\left(\mathcal{A}(K_c)M_c\right)^{\abs{\eta}}\abs{\eta}!, \;\; \forall\eta\in X^\ast,
\label{eq:local-unity-feedback-growth-bound}
\endeq
for some $K>0$ and
\begdi
\mathcal{A}(K_c)=\frac{1}{1-mK_c \ln\left(1+\frac{1}{mK_c}\right)}.
\enddi
Therefore, $c^{-1}\in\allseriesmLC$.
Furthermore, no geometric growth constant smaller than $\mathcal{A}(K_c)M_c$
can satisfy \rref{eq:local-unity-feedback-growth-bound}, so the radius
of convergence for this operator is $1/\mathcal{A}(K_c)M_c(m+1)$.
\endth

\begpr
It was shown in \cite[Corollary 2]{Thitsa-Gray_12} that the generating series
for the unity
feedback system $c@\delta$
has exactly the properties described above, and
therefore, so does $(-c)@\delta$. The present theorem
is thus proved by showing that $c^{-1}=(-c)@\delta$.
Recall that in the proof of  Theorem~\ref{th:allseriesdeltam-is-group} it was
shown in general that $c^{-1}=(-c)\modcomp c^{-1}$.
But it is also known that
$(-c)@\delta$ satisfies the fixed point equation
$(-c)@\delta=(-c)\modcomp((-c)@\delta)$ \cite{Gray-Li_05}.
Therefore, since $e\mapsto (-c)\modcomp e$ is a contraction on a complete
ultrametric space, the identity in question must
hold.
\endpr

It is worth noting that
the growth constants determined in Theorem~\ref{th:ROC-for-cinv_LC} must hold for
every series $c$ with growth constants $K_c,M_c$.
Thus, it tends to be conservative for specific series in
this class (see \cite{Thitsa-Gray_12} for further discussion on this topic).
A similar approach yields the global counterpart of this theorem.

\begth \label{th:ROC-for-cinv_GC}
For any $c\in\allseriesmGC$ with
growth constants $K_c,M_c>0$ it follows that
\begeq
|(c^{-1},\eta)|\leq K\left(\mathcal{B}(K_c)M_c\right)^{\abs{\eta}}\abs{\eta}!, \;\; \forall\eta\in X^\ast,
\label{eq:global-unity-feedback-growth-bound}
\endeq
for some $K>0$ and
\begdi
\mathcal{B}(K_c)=\frac{1}{\ln\left(1+\frac{1}{mK_c}\right)}.
\enddi
Therefore, $c^{-1}\in\allseriesmLC$.
Furthermore, no geometric growth constant smaller than $\mathcal{B}(K_c)M_c$
can satisfy \rref{eq:global-unity-feedback-growth-bound}, so the radius
of convergence for this operator is $1/\mathcal{B}(K_c)M_c(m+1)$.
\endth

It is known that feedback does not in general preserve global convergence (see \cite{Gray-etal_09} for a
specific example). Thus, there is no reason to expect that the composition inverse will do so either.

\section{Feedback Product}

The goal of this section is to derive an explicit formula for the multivariable feedback product
$c@d$ using the Fa\`{a} di Bruno Hopf algebra described in the previous section.
Given two Fliess operators $F_c$ and $F_d$ which are linear time-invariant
systems with $\ell_c\times m_c$ transfer function $G_c$ and
$\ell_d\times m_d$ transfer function $G_d$, respectively, the closed-loop
transfer function is clearly
\begeq
G_c(I-G_dG_c)^{-1}=G_c\sum_{k=0}^\infty (G_dG_c)^k, \label{eq:MIMO-CL-transfer-function}
\endeq
where necessarily $\ell_c=m_d$ and $\ell_d=m_c$. There is no a priori
requirement that the systems be {\em square}, that is, $m_c=\ell_c$ or
$m_d=\ell_d$. But to handle the most general case here, the series composition products
introduced in Sections~\ref{sec:preliminaries} and \ref{sec:FdB-Hopf-algebra}
have to be generalized to accommodate two alphabets, $X_c=\{x_0,x_1,\ldots,x_{m_c}\}$
and $X_d=\{\tilde{x}_0,\tilde{x}_1,\ldots,\tilde{x}_{m_d}\}$. This offers no serious technical
issues as described in \cite[Example 3.5]{Gray-Thitsa_IJC12}, just a bit more bookkeeping.
The inverse is computed easily in this special case because all the underlying series are
rational.
The next theorem gives the nonlinear generalization of \rref{eq:MIMO-CL-transfer-function}.

\begth
\label{th:feedback-product-formula}
For any $c\in\allseriesXcmd$ and $d\in\allseriesXdmc$, it follows that
\begeq
c@d=c\modcomp(-d\circ c)^{-1}=c\circ(\delta-d\circ c)^{-1}. \label{eq:catd-formula} \\[0.1in]
\endeq
\endth

\begpr
The proof is not significantly different from the SISO case presented in \cite{Gray-Duffaut_Espinosa_SCL11}.
Since it is short, it is presented here for completeness.
Clearly the function $v$ in Figure~\ref{fig:feedback-with-v} must
satisfy the identity
\begdi
v=u+F_{d\circ c}[v].
\enddi
Therefore,
\begdi
\left(I+F_{-d\circ c}\right)[v]=u,
\enddi
where in the notation of Section~\ref{sec:FdB-Hopf-algebra}, the operator on the left-hand side
is an element of $\Fliessdelta$ with $m=m_c=\ell_d$.
Applying the composition inverse $\left(I+F_{(-d\circ c)^{-1}}\right)$
on the left gives
\begdi
v=\left(I+F_{(-d\circ c)^{-1}}\right)[u],
\enddi
and thus,
\begdi
F_{c@d}[u]=F_c[v]=F_{c\modcomp(-d\circ c)^{-1}}[u]
\enddi
as desired. The second identity in the theorem
is just a formal way of expressing the first identity
since $c\modcomp(-d\circ c)^{-1}=c\circ (\delta+(-d\circ c)^{-1})$ and by definition
$(\delta-d\circ c)^{-1}=\delta+(-d\circ c)^{-1}$.
\endpr

As noted earlier, \rref{eq:catd-formula} also makes sense when either
$c=\delta$ or $d=\delta$, namely,
$\delta@d=(\delta-d)^{-1}=\delta+(-d)^{-1}$ and
$c@\delta=c\circ(\delta-c)^{-1}=(-c)^{-1}$. In addition, it was shown in \cite[Theorem 4.3]{Gray-Li_05}
that $c@d$ satisfies the fixed point equation $c@d=c\circ (\delta+d\circ(c@d))$. So if $c$ is a linear series then
\begin{align*}
c@d&=c+c\circ d\circ(c@d) \\
(\delta-c\circ d)\circ (c@d)&=c\\
c@d&=(\delta-c\circ d)^{-1}\circ c.
\end{align*}
But in general, even in the SISO case, $c@d\neq (\delta-c\circ d)^{-1}\circ c$.

Next it is shown that feedback preserves local convergence. But the following
preliminary result is needed first.

\begth \label{th:FdB-Hopf-Fliess-subgroup}
The triple $(\allseriesdeltamLC,\circ,\delta)$ is a subgroup of
$(\allseriesdeltam,\circ,\delta)$.
\endth
\begpr
The set of series $\allseriesdeltamLC$ is
closed under composition since the set $\allseriesmLC$ is
closed under addition and modified composition \cite{Gray-Li_05,Li_04}.
In light of Theorem~\ref{th:ROC-for-cinv_LC},
$\allseriesmLC$ is also closed under inversion. Hence, the theorem is proved.

\endpr

\begth
If $c\in\allseriesXcmdLC$ and $d\in\allseriesXdmcLC$  then
$c@d\in\allseriesXcmdLC$.
\endth

\begpr
Since the composition product, the modified composition product, and the
composition inverse all preserve local convergence, the claim follows directly from
Theorem~\ref{th:feedback-product-formula}.

\endpr

\begex
{\rm Consider the differential axle shown in Figure \ref{fig:Diff_axle}.
This device has zero mass and moves in the plane with independent
angular velocities $u_r$ and $u_l$ corresponding to the right and left
wheels, respectively.
\begin{figure}[t]
\begce
\includegraphics[width=8cm]{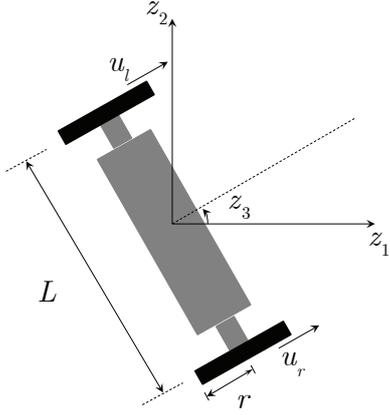}
\endce
\caption{Differential axle}
\label{fig:Diff_axle}
\end{figure}
The dynamics of this system are
\begin{align*}
\dot{z_1} & =  \frac{r}{2} (u_l+u_r) \cos(z_3)\\
\dot{z_2} & =  \frac{r}{2} (u_l+u_r) \sin(z_3)\\
\dot{z_3} & =  \frac{r}{L} (u_r-u_l).
\end{align*}
In particular, if $u_l=u_r>0$ then the axle moves forward in the direction the
wheels are pointing, and if $u_l=-u_r>0$ the axle rotates clockwise because the
wheels are turning in opposite directions. For simplicity, define $u_1=
\frac{1}{2}(u_l+u_r)$ and $u_2=(u_r-u_l)$, and let
$L=r=1$. Choosing outputs $y_i= z_i$, $i=1,2$, the corresponding
two-input, two-output state space realization is
\begin{subequations} \label{eq:simplified_diff-system}
\begin{align}
\left(\begin{array}{c} \dot{z_1} \\ \dot{z_2} \\ \dot{z_3} \end{array}\right) &
= \left(\begin{array}{c} \cos(z_3) \\ \sin(z_3) \\ 0 \end{array}\right) u_1 +
\left(\begin{array}{c} 0 \\ 0 \\ 1 \end{array}\right) u_2 \\
\left(\begin{array}{c} y_1 \\ y_2 \end{array}\right) & =
\left(\begin{array}{c} z_1 \\ z_2 \end{array}\right).
\end{align}
\end{subequations}
Its generating series, $c$, can be computed directly from \rref{eq:c-equals-Lgh}
using the vector fields and output function given in
\rref{eq:simplified_diff-system}. With the help of the
Mathematica software package NCAlgebra \cite{NCA_12},
this calculation gives
\begdi
c=\left(\begin{array}{c}
           x_1 - x_1x_2^2 + x_1x_2^4 -  x_1x_2^6 +  x_1x_2^8
-x_1x_2^{10}+\cdots\\
           x_1x_2  - x_1x_2^3  + x_1x_2^5 - x_1x_2^7 +
x_1x_2^9-x_1x_2^{11}+\cdots
\end{array}\right)
\enddi
when $z_1(0)=z_2(0)=z_3(0)=0$. This series is clearly in $\re_{GC}^2\langle
\langle X \rangle\rangle$ with $X=\{x_0,x_1,x_2\}$
and growth constants $K_c=1$ and $M_c=1$.

Consider now the problem of steering the differential axle around a circle. For this
purpose, a two channel proportional-integral controller is used in the feedback
path so that one obtains a closed-loop system as shown in
Figure~\ref{fig:feedback-with-v}. The dynamics of the
controller are
\begin{subequations} \label{subeq:Differentialaxle_output}
\begin{align}
\left(\begin{array}{c} \dot{z_4} \\ \dot{z_5} \end{array}\right) &
= \left(\begin{array}{c} 1 \\ 0  \end{array}\right) \tilde{u}_1
+ \left(\begin{array}{c} 0 \\ 1  \end{array}\right) \tilde{u}_2 \\
\left(\begin{array}{c} \tilde{y}_1 \\ \tilde{y}_2 \end{array}\right) & =
\left(\begin{array}{c} k_1z_4 \\ k_2z_5 \end{array}\right).
\end{align}
\end{subequations}
For gains $k_1=2$ and $k_2=10$, the corresponding generating
series is
\begdi
d=\left(\begin{array}{c}
           4+ 2 x_1 \\
           20 + 10 x_2
\end{array}\right)
\enddi
when $z_4(0)=z_5(0)=2$. Here $d\in \re_{GC}^2\langle \langle X
\rangle\rangle$ with growth constants $K_d= 20$ and
$M_d= 0.5$.

From Theorem \ref{th:feedback-product-formula}, the feedback product of
$c$ and $d$ is given by the series
\begin{align*}
\lefteqn{(c @ d)_1=} \\
& 4 x_0 + x_1 - 1,592 x_0^3 + 2 x_0^2 x_1 - 80 x_0^2 x_2 - 80 x_0 x_2 x_0 - \\
& 4 x_0 x_2^2 - 400 x_1 x_0^2 - 20 x_1 x_0 x_2 - 20 x_1 x_0 x_2 - \\
& 20 x_1 x_2 x_0 -x_1 x_2^2 + 617,616 x_0^5 -2,396 x_0 ^4 x_1 + \\
& 31,360 x_0^4 x_2 - 1,600 x_0^3 x_1 x_0 - 80 x_0^3 x_1 x_2 + 31,360 x_0^3 x_2 x_0 - \\
& 80 x_0^3 x_2 x_1  + 1,584 x_0^3 x_2^2 -1,600 x_0^2 x_1 x_0^2- 80 x_0^2 x_1 x_0 x_2 - \\
&80 x_0^2 x_1 x_2 x_0 - 4 x_0^2 x_1 x_2^2
 + 31,520 x_0^2 x_2 x_0^2
- 80 x_0^2 x_2 x_0 x_1 + \\
& 1,592 x_0^2 x_2 x_0 x_2
 - 40 x_0^2 x_2 x_1 x_0
 - 2 x_0^2 x_2 x_1 x_2 + \\
& 1,592 x_0^2 x_2^2 x_0
 -2 x_0^2 x_2^2 x_1 + 80 x_0^2 x_2^3
 + 31,520 x_0 x_2 x_0^3- \\
& 80 x_0 x_2 x_0^2 x_1 + 1,592 x_0 x_2 x_0^2 x_2
 - 40 x_0 x_2 x_0 x_1 x_0 - \\
& 2 x_0 x_2 x_0 x_1 x_2
 + 1,592 x_0 x_2 x_0 x_2 x_0
 -2 x_0 x_2 x_0 x_2 x_1+ \\
& 80 x_0 x_2 x_0 x_2^2
 + 1,592 x_0 x_2^2 x_0^2
  - 2 x_0 x_2^2 x_0 x_1 +
160,000 x_1 x_0^4 + \\
&\cdots
\end{align*}
and
\begin{align*}
\lefteqn{(c @ d)_2=} \\
& 80 x_0^2 + 4x_0x_2 +20 x_1 x_0 + x_1 x_2 - 31,520 x_0^4
 + 80 x_0^3 x_1 - \\
& 1,592 x_0^3 x_2 + 40 x_0^2 x_1 x_0
 + 2 x_0^2 x_1 x_2
 - 1,592 x_0^2 x_2 x_0 + \\
& 2 x_0^2 x_2 x_1
 - 80 x_0^2 x_2^2 - 1,592 x_0 x_2 x_0^2
 +2 x_0x_2 x_0 x_1 -\\
&80 x_0 x_2 x_0 x_2 - 80 x_0 x_2^2 x_0 - 4 x_0 x_2^3
 - 8,000 x_1 x_0^3 - \\
& 400 x_1 x_0^2 x_2 - 400 x_1 x_0 x_2 x_0
-20 x_1x_0
x_2^2
 -400 x_1 x_2 x_0^2 - \\
& 20 x_1 x_2 x_0 x_2
 - 20 x_1 x_2^2 x_0
 - x_1 x_2^3 + 3,200 x_0^5 + 160 x_0^4 x_2+ \\
&  800 x_0^3 x_1 x_0 + 40 x_0^3 x_1 x_2 + 800 x_1 x_0^4
 + 40 x_1 x_0^3 x_2+ \\
&  200 x_1x_0^2 x_1 x_0 + 10 x_1 x_0^2 x_1 x_2
 + 11,841,600 x_0^6 +\cdots.
\end{align*}%
\begin{figure}[t]
\hspace*{-0.1in}\includegraphics[width=9.2cm]{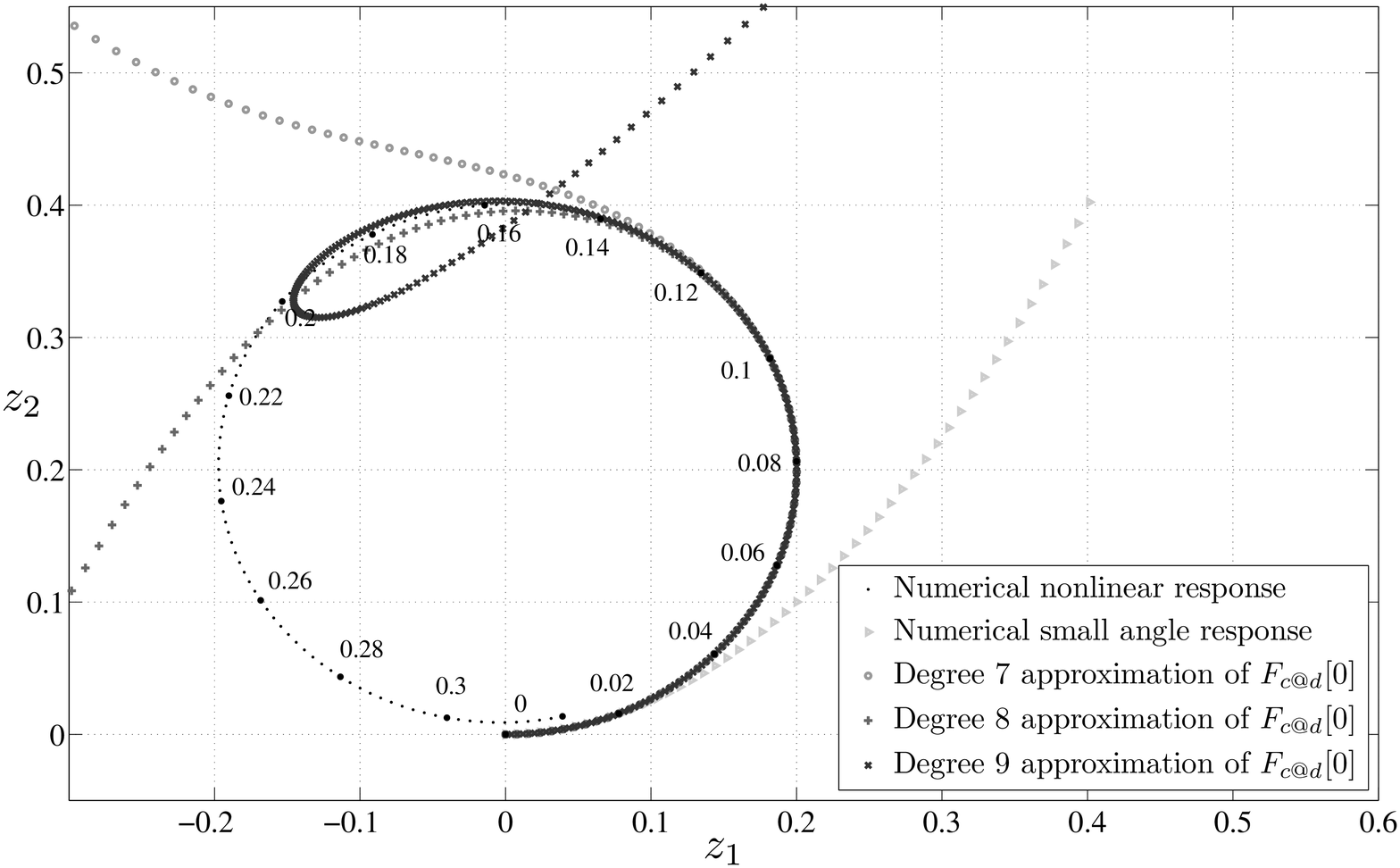}
\caption{Various estimates of the natural response of the closed-loop
differential axle system}
\label{fig:Diff_axle_output}
\end{figure}%
The outputs of the closed-loop system for a zero reference input are then
\begin{align*}
y_{1}(t) & =  F_{(c@ d)_{1}}[0](t) \\
&=  4t\! - \frac{796 t^3}{3}\! +
\frac{25,734t^5}{5}\! - \frac{4000 t^6}{9}\!
 - \frac{13,528,798 t^7}{315}+ \\
 &\hspace*{0.18in} \frac{1,653,800 t^8}{63}\! +
\frac{594,150,001 t^9}{5670}\! + \cdots
\end{align*}
and
\begin{align*}
y_{2}(t) & = F_{(c@ d)_{2}}[0](t) \\
&= 40 t^2\! - \frac{3,940
t^4}{3}\! + \frac{80 t^5}{3}\! +
\frac{49,340 t^6}{3}\! - \frac{254,560 t^7}{63} -\\
&\hspace*{0.18in} \frac{5,496,593 t^8}{63}\! + \frac{25,525,060 t^9}{189}\! +\cdots
\end{align*}
Various estimates of the natural response of the closed-loop differential axle
system are shown in Figure \ref{fig:Diff_axle_output}.
The tick marks along the circle indicate time.
Specifically, the
numerically computed nonlinear closed-loop response of
\rref{eq:simplified_diff-system}-\rref{subeq:Differentialaxle_output}
is compared against Fliess operator responses whose generating series are
computed from the feedback product truncated to degrees
$7$, $8$ and $9$. Also shown in Figure~\ref{fig:Diff_axle_output} is
the response of the small angle approximation system
\begin{subequations}
\begin{align*}
\left(\begin{array}{c} \dot{z_1} \\ \dot{z_2} \\ \dot{z_3} \end{array}\right) &
= \left(\begin{array}{c} 1 \\ z_3 \\ 0 \end{array}\right) u_1 +
\left(\begin{array}{c} 0 \\ 0 \\ 1 \end{array}\right) u_2 \\
\left(\begin{array}{c} y_1 \\ y_2 \end{array}\right) & =
\left(\begin{array}{c} z_1 \\ z_2 \end{array}\right)
\end{align*}
\end{subequations}
steered by the same proportional-integral controller.
It is evident that this system
underestimates the correct position of the differential axle
almost immediately. On the other hand, the Fliess operator approximations
clearly improve as additional terms are added to the approximation.

One way to get some insight into the convergence characteristics of $F_{c@d}[0]$
is to empirically estimate the geometric growth constants for the natural response
portion of the series $c@d$, i.e., the series $(c@d)_N:=\sum_{k\geq 0} (c@d,x_0^k)x_0^k$,
by plotting $\ln((c@d)_N/\abs{\eta}!)$ versus $\abs{\eta}$ for the local case
and $\ln((c@d)_N)$ versus $\abs{\eta}$ for the global case as shown in
Figure~\ref{fig:yN-growth-constant-plots}. To improve the quality of the estimates,
coefficients above order nine were computed using \rref{eq:c-equals-Lgh}.
In each case, the corresponding growth constants
can be estimated by linearly fitting the data
(see \cite{Gray-Duffaut_Espinosa_FdB14} for more discussion concerning this methodology).
The parameter $R^2$ is the square of Pearson's correlation coefficient, so the closer this
statistic is to unity, the
better the linear fit. In this case, the data appears to match better the global growth rate with
$M_{(c@d)_N}=\exp(3.1157)=22.549$. But as will be discussed shortly, the series can fall
somewhere {\em in between} being locally convergent and globally convergent
as defined by \rref{eq:local-convergence-growth-bound} and
\rref{eq:global-convergence-growth-bound}, respectively.
\begin{figure}[t]
\hspace*{-0.1in}\includegraphics[width=9.2cm]{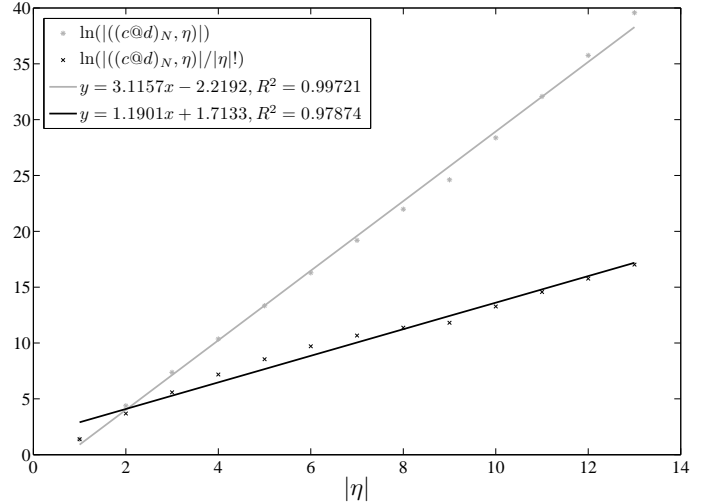}
\caption{Linear fits of the coefficients of $(c@d)_N$ on a logarithmic scale}
\label{fig:yN-growth-constant-plots}
\end{figure}
There is also
the option of constructing a piecewise analytic approximation of the response
using a sequence of closed-loop generating series computed by brute force or
via analytic extension \cite{Moser_96}. This approach has the additional
advantage that lower order approximations of each piece are likely to suffice. For example,
it appears here that two
degree nine approximations joined at the midpoint of the path could easily traverse the entire
circle.

An alternative method to exploring the nature of the convergence of $F_{c@d}[0]$ is to use
Theorems~\ref{th:ROC-for-cinv_LC} or \ref{th:ROC-for-cinv_GC} in conjunction with what
is known at present about the convergence of interconnected Fliess operators
as reported in \cite{Thitsa-Gray_12}.
First observe that
$(-d\circ c, \eta)=\pm{2 \choose 10}$ when $\eta\neq \emptyset$, and therefore,  $-d\circ c$ is globally convergent
with
$K_{d\circ c}=20$ and $M_{d\circ c}=1$.
Applying Theorem~\ref{th:ROC-for-cinv_GC} provides
an upper bound on the
local geometric growth constant of $e:=(-d\circ c)^{-1}$, specifically,
\begin{align*}
M_{e}&=\frac{M_{d\circ c}}{\ln\left(1+\frac{1}{2K_{d\circ c}}\right)}=40.50.
\end{align*}
Repeating the empirical method used above for $e$ gives
the data shown in Figure~\ref{fig:e-growth-constant-plots}.
\begin{figure}[t]
\hspace*{-0.1in}\includegraphics[width=9.2cm]{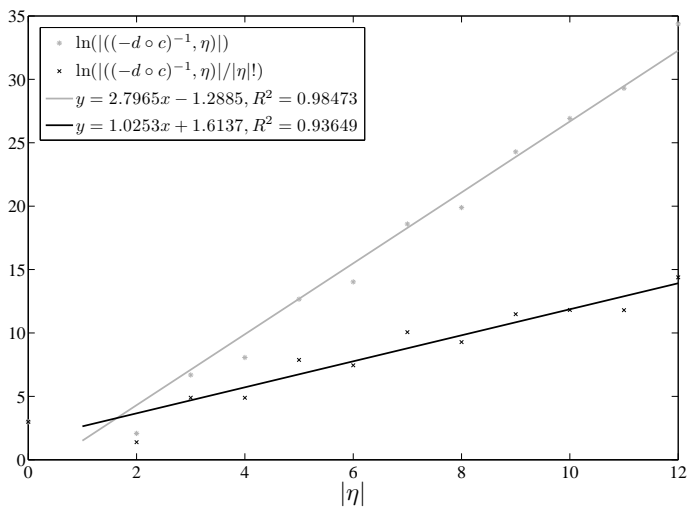}
\caption{Linear fits of the coefficients of $e=(-d\circ c)^{-1}$ on a logarithmic scale}
\label{fig:e-growth-constant-plots}
\end{figure}
It indicates that $e$ is more globally convergent in nature than locally
convergent, but if it
were assumed to be locally convergent,
the corresponding geometric growth constant would be $\exp(1.0253)=2.7879<40.50$.
On the other hand, if it were taken to be globally convergent then (since
$F_{c\tilde{\circ} e}[0]=F_{c\circ e}[0]$) it
follows that $(c@d)_N=c\circ e$ is the composition of two globally
convergence series. As discussed in \cite[p.~2800]{Thitsa-Gray_12}, the
resulting series can lie strictly in between locally and globally
convergent.
But independent of this fact, it is still known in this instance that the series
$F_{c\circ e}[0]$ will
converge over any finite interval \cite[Theorem 9]{Thitsa-Gray_12}.
}
\endex

\section{Conclusions}

The main thrust of this paper was to provide the
full multivariable extension of a theory to explicitly compute the
generating series of a feedback interconnection of
two systems represented as Fliess operators.
This was largely facilitated by
utilizing
a new type of grading for the underlying Hopf algebra.
This grading also provided a fully recursive
algorithm to compute the antipode of the algebra
and thus, the corresponding
feedback product can be computed much more efficiently.
Finally, an improved convergence analysis of
the antipode operation was presented, one that gives the radius of convergence
for this operation.

\vspace*{0.1in}

\noindent
{\bf Acknowledgements}

\vspace*{0.1in}

The first author was supported by grant SEV-2011-0087 from the
Severo Ochoa Excellence Program at the Instituto de Ciencias Matem\'{a}ticas
in Madrid, Spain. The third author was supported by Ram\'on y Cajal research
grant RYC-2010-06995 from the Spanish government. The authors thank the reviewers for their
suggestions to improve the presentation.


\begin{thebibliography}{99}

\bibitem{Duffaut-Espinosa_09}
 L.~A.~Duffaut Espinosa, Interconnections of Nonlinear Systems Driven by $L_2$-It\^{o} Stochastic
Processes, Doctoral dissertation, Old Dominion University, 2009.

\bibitem{Duffaut-Espinosa-et-al_CDC09}
L.~A.~Duffaut Espinosa, W.~S.~Gray, and O.~R.~Gonz\'{a}lez,
On Fliess operators driven by $L_2$-It\^o random processes,
{\em Proc.~48th IEEE Conf. on Decision and Control},
Shanghai, China, 2009, pp.~7478--7484.

\bibitem{Einziger_10}
H. Einziger,
Incidence Hopf Algebras: Antipodes, Forest Formulas, and Noncrossing Partitions,
Doctoral dissertation, The George Washington University, 2010.

\bibitem{Ferfera_79}
A.~Ferfera, Combinatoire du Mono\"{i}de Libre Appliqu\'{e}e
\`{a} la Composition et aux Variations de Certaines Fonctionnelles
Issues de la Th\'{e}orie des Syst\`{e}mes, Doctoral dissertation,
University of Bordeaux I, 1979.

\bibitem{Ferfera_80}
\sameau, Combinatoire du mono\"{i}de libre et composition de
certains syst\`{e}mes non lin\'{e}aires, {\em Ast\'{e}risque},
75--76 (1980) 87--93.

\bibitem{Figueroa-Gracia-Bondia_05}
H.~Figueroa and J.~M.~Gracia-Bond\'{ii}a,
Combinatorial Hopf algebras in quantum field theory I,
{\em Rev.~Math.~Phys.}, 17 (2005) 881--976.

\bibitem{Fliess_81}
M.~Fliess, Fonctionnelles causales non lin\'{e}aires et ind\'{e}termin\'{e}es non
commutatives, {\em Bull.~Soc.~Math.~France}, 109 (1981) 3--40.

\bibitem{Fliess_83}
\sameau, R\'{e}alisation locale des syst\`{e}mes non lin\'{e}aires, alg\`{e}bres de
Lie filtr\'{e}es transitives et s\'{e}ries g\'{e}n\'{e}ratrices non commutatives,
{\em Invent.~Math.}, 71 (1983) 521--537.

\bibitem{Foissy_13}
L.~Foissy,
The Hopf algebra of Fliess operators and its dual pre-Lie algebra,
{\tt http://arxiv.org/abs/1304.1726v3}, 2014.

\bibitem{Gray_MTNS14}
W. S. Gray, Affine feedback transformation group for nonlinear SISO systems,
{\em Proc.\ 21st Inter.~Symp.~Mathematical Theory
of Networks and Systems}, Groningen, The Netherlands, 2014, pp.~297--302.

\bibitem{Gray-Duffaut_Espinosa_SCL11}
W. S. Gray and L. A. Duffaut Espinosa, A Fa\`{a} di Bruno Hopf algebra for a group of
Fliess operators with applications to feedback, {\em Systems Control Lett.}, 60 (2011) 441--449.

\bibitem{Gray-Duffaut_Espinosa_FdB14}
\sameau,
A Fa\`{a} di Bruno Hopf algebra for analytic nonlinear feedback control systems,
in {\em Fa\`{a} di Bruno Hopf Algebras, Dyson-Schwinger
Equations, and Lie-Butcher Series}, K. Ebrahimi-Fard and F. Fauvet, Eds.,
IRMA Lect.~Math.~Theor.~Phys.,
Eur.~Math.~Soc.,
Strasbourg, France, to appear.

\bibitem{Gray-et-al_MTNS14}
W.~S.~Gray, L.~A.~Duffaut Espinosa, and K.~Ebrahimi-Fard,
Recursive algorithm for the antipode in the SISO feedback product,
{\em Proc.\ 21st Inter.~Symp.~Mathematical Theory
of Networks and Systems}, Groningen, The Netherlands, 2014, pp.~1088--1093.

\bibitem{Gray-et-al_AUT}
W. S. Gray, L. A. Duffaut Espinosa, and M. Thitsa,
Left inversion of analytic nonlinear SISO systems via formal power series methods,
{\em Automatica}, 50 (2014) 2381-–2388.

\bibitem{Gray-etal_09}
W. S. Gray, H. Herencia-Zapana, L. A. Duffaut Espinosa,
and O. R. Gonz\'{a}lez,
Bilinear system interconnections and
generating series of weighted Petri nets,
{\em Systems Control Lett.}, 58 (2009) 841--848.

\bibitem{Gray-Li_05}
W.~S.~Gray and Y.~Li, Generating series for interconnected analytic
nonlinear systems, {\em SIAM J.~Control Optim.}, 44 (2005) 646--672.

\bibitem{Gray-Thitsa_IJC12}
W.~S.~Gray and M.~Thitsa,
A unified approach to generating series for mixed cascades of
analytic nonlinear input-output systems,
{\em Inter. J. Control}, 85 (2012) 1737--1754.

\bibitem{Gray-Wang_SCL02}
W.~S.~Gray and Y.~Wang, Fliess operators on $L_p$ spaces:
convergence and continuity, {\em Systems Control Lett.}, 46 (2002) 67--74.

\bibitem{Gray-Wang_08}
\sameau,
Formal Fliess operators with applications to feedback interconnections,
{\em Proc.~18th Inter.~Symp.~Mathematical Theory of
Networks and Systems}, Blacksburg, Virginia, 2008.

\bibitem{Haiman-Schmitt_89}
M.~Haiman and W.~Schmitt,
Incidence algebra antipodes and Lagrange inversion
in one and several variables,
{\em J. Combin. Theory Ser. A},
50 (1989) 172--185.

\bibitem{Hochschild_81}
\markblue{
G. P. Hochschild, {\em Basic Theory of Algebraic Groups and Lie Algebras}, Springer-Verlag, New York, 1981.}

\bibitem{Isidori_95}
A.~Isidori, {\em Nonlinear Control Systems}, 3rd Ed., Springer-Verlag, London, 1995.

\bibitem{Li_04}
Y.~Li, Generating Series of Interconnected Nonlinear Systems and the
Formal {L}aplace-{B}orel Transform,
Doctoral dissertation, Old Dominion University, 2004.

\bibitem{Moser_96}
A. Moser,
Extending the domain of definition of functional series
for nonlinear systems,
{\em Automatica}, 32 (1996) 1233--1234.

\bibitem{NCA_12}
The NCAlgebra Suite, Version 4.0, currently available at {\tt math.} {\tt ucsd.edu/$\sim$ncalg},
2012.

\bibitem{Reutenauer_93}
C.~Reutenauer,
{\em Free Lie Algebras},
Oxford University Press, New York, 1993.

\bibitem{Sweedler_69}
M.~E.~Sweedler,
{\em Hopf Algebras}, W.~A.~Benjamin, Inc., New York, 1969.

\bibitem{Thitsa-Gray_12}
M. Thitsa and W. S. Gray,
On the radius of convergence of interconnected
analytic nonlinear input-output systems,
{\em SIAM J.~Control Optim.}, 50 (2012) 2786--2813.

\end{thebibliography}
\end{document}